\documentclass[pdflatex,sn-mathphys-num]{sn-jnl}

\usepackage{graphicx}%
\usepackage{multirow}%
\usepackage{amsmath,amssymb,amsfonts}%
\usepackage{amsthm}%
\usepackage{mathrsfs}%
\usepackage[title]{appendix}%
\usepackage{xcolor}%
\usepackage{textcomp}%
\usepackage{manyfoot}%
\usepackage{booktabs}%
\usepackage{algorithm}%
\usepackage{algorithmicx}%
\usepackage{algpseudocode}%
\usepackage{listings}%

\usepackage{diagbox}
\usepackage{makecell}
\usepackage{tabularx}

\newcolumntype{C}{>{\centering\arraybackslash}X}

\begin{document}

\title[Article Title]{A piecewise neural network method for solving large interval solution to initial value problem of ordinary differential equations}


\author[1]{\fnm{Dongpeng} \sur{Han}}\email{handongpeng@stu.shmtu.edu.cn}

\author*[2]{\fnm{Chaolu} \sur{Temuer}}\email{tmchaolu@shmtu.edu.cn}
\equalcont{These authors contributed equally to this work.}

\affil[1,2]{\orgdiv{Mathematics Department}, \orgname{Shanghai Maritime University}, \orgaddress{\street{1550 Haigang Avenue}, \city{Shanghai}, \postcode{201306}, \state{Shanghai}, \country{China}}}


\abstract{Various traditional numerical methods for solving initial value problems of differential equations often produce local solutions near the initial value point, despite the problems having larger interval solutions. Even current popular neural network algorithms or deep learning methods cannot guarantee yielding large interval solutions for these problems. In this paper, we propose a piecewise neural network approach to obtain a large interval numerical solution for initial value problems of differential equations. In this method, we first divide the solution interval, on which the initial problem is to be solved, into several smaller intervals. Neural networks with a unified structure are then employed on each sub-interval to solve the related sub-problems. By assembling these neural network solutions, a piecewise expression of the large interval solution to the problem is constructed, referred to as the piecewise neural network solution. The continuous differentiability of the solution over the entire interval, except for finite points, is proven through theoretical analysis and employing a parameter transfer technique. Additionally, a parameter transfer and multiple rounds of pre-training technique are utilized to enhance the accuracy of the approximation solution. Compared with existing neural network algorithms, this method does not increase the network size and training data scale for training the network on each sub-domain. Finally, several numerical experiments are presented to demonstrate the efficiency of the proposed algorithm.}


\keywords{Piecewise neural network, Initial value  problem, Ordinary differential equations, Extension of local solution, Parameter transfer}



\maketitle

\section{Introduction}
As we know, constrained by convergence requirements, various traditional numerical solution methods, such as Picard approximation and Runge-Kutta methods, for initial value problems (IVP) of differential equations or dynamic systems (ODEs), are typically suitable only for local solutions within a neighborhood interval of the initial value point, even if the problem possesses larger interval (global) solutions. Consequently, extending a local solution to a larger interval remains a fundamental challenge in the field of numerical solutions to ODE problems, which has yet to be fully addressed \cite{ODE, Maria}.

Currently, the artificial neural network (ANN) algorithm, as a novel numerical solution method for solving differential equation problems, is rapidly evolving and attracting attention and exploration from many researchers \cite{NUrODE}. Although the ANN algorithm, like traditional methods, performs well mainly in the neighborhood of the initial value point and lacks inherent learning capability on large intervals, it offers numerous advantages in solving IVPs of ODEs. For instance, the neural network-based solution of a differential equation is differentiable and presented in closed analytic form, suitable for subsequent calculations. It is independent of discrete schemes and the shape of the variable domain. Therefore, it is anticipated that the aforementioned solution extension problem can be effectively addressed through the enhancement and efficient utilization of the ANN algorithm. In this paper, we propose a piecewise ANN (PWNN) method to tackle this problem.

Before describing our algorithm, let us briefly review the development profile of ANN algorithms for solving differential equation problems.

The theoretical foundation supporting the use of ANN algorithms for solving differential equation problems lies in the general approximation theorems provided by Hornik, Womik, and Li X et al. \cite{hornik_multilayer_1989, z49, z50}, which theoretically establish that any continuous function can be approximated by a neural network. Subsequently, a widely utilized ANN algorithm for solving differential equations with initial boundary value conditions was introduced by Lagaris et al. \cite{lagaris_artificial_1998}. However, for problems involving more complex initial boundary values, this method faces difficulties in constructing the necessary trial solutions for the underlying problems.

Another significant development in ANN algorithms for solving differential equation problems is the implementation of automatic differentiation technology \cite{paszke_pytorch_2019, z46, z47, z48}, which has led to the proliferation of ANN methods for solving differential equation problems and their broader application. Building upon this technique, Justin Sirignano and Konstantinos Spiliopoulos proposed the Deep Galerkin method (DGM), a deep neural network method for solving (partial) differential equation problems in higher dimensions \cite{sirignano_dgm_2018}. Cosmin Anitescu et al. introduced an ANN method utilizing adaptive configuration strategies to enhance method robustness and reduce computational costs in solving boundary value problems of differential equations \cite{anitescu_artificial_2019}.

A landmark advancement in ANN algorithms for solving differential equation problems was the introduction of physical information neural networks (PINN) by Raissi et al. \cite{raissi_physics-informed_2019}. In this method, the initial boundary conditions of the differential equation are incorporated into the loss function, allowing the ANN to directly express an approximate solution to the differential equation. This eliminates the need to construct a trial solution according to equations and initial boundary value conditions, as required by Lagaris's method. Subsequently, numerous ANN algorithms based on PINN have been developed rapidly. For example, Lei Yuan and Yi-Qing Ni et al. introduced the Auxiliary PINN (A-PINN), capable of bypassing limitations in integral discretization and solving forward and inverse problems of nonlinear integral differential equations \cite{yuan_-pinn_2022}. Pao-Hsiung Chiu et al. proposed novel PINN methods for coupling neighboring support points and their derivative terms obtained by automatic differentiation \cite{chiu_can-pinn_2022}. Yao Huang et al. combined PINN with the homotopy continuation method, proposing Homotopy PINN (HomPINN) for solving multiple solutions of nonlinear elliptic differential equations, overcoming the limitation of PINN in finding only the flattest solution in most cases \cite{huang_hompinns_2022}. Fang, Yin et al. utilized PINN to address a range of femtosecond optical soliton solutions pertaining to the high-order nonlinear Schr${\ddot{\textrm{o}}}$dinger equation \cite{fang_data-driven_2021}. Based on the PINN method, Bai, Yuexing et al. proposed an enhanced version of PINN called IPINN, introducing localized adaptive activation functions to improve performance, successfully applying the method to solve several differential equation models in finance \cite{bai_application_2022}. Meng et al. introduced a modified PINN method called PPINN, dividing a long-period problem into a series of short-period ones to accelerate the training of ANN algorithms \cite{z80}.

Moreover, there are numerous other analogues of PINN-based studies and various types of ANN methods. For instance, convolutional ANN methods \cite{z57,z73, z93} and theoretically guided ANN methods \cite{z85, z86} have also been employed to explore new solution methods for (partial) differential equations. Particularly noteworthy is the work of Run-Fa Zhang and Sudao Bilige, who proposed bilinear neural networks, the first attempt to obtain analytical solutions to nonlinear partial differential equations using the ANN method \cite{zhang_bilinear_2019,zhang2}.

Importantly, a general observation from these literature is that ANN algorithms for solving IVPs of differential equations perform well in relatively small domains near the initial value point but sometimes exhibit poor convergence over large intervals. To obtain a large interval solution, it is often necessary to enhance the training data and expand the network size (depth and width). However, doing so not only reduces computational efficiency but also makes it challenging to guarantee obtaining a large interval solution to the problem. Therefore, alternative strategies are being explored to enhance the efficiency of ANN algorithms. For example, various adaptive activation functions \cite{bai_application_2022, chien-cheng_yu_adaptive_2002, dushkoff_adaptive_2016, qian_adaptive_2018} and adaptive weights \cite{wang_when_2022,xiang_self-adaptive_2022} have been proposed. The key feature of these studies is the introduction of hyperparameters into traditional activation functions to adjust ANN convergence. Studies by Jaftap et al. have shown that these methods are primarily effective in the early stages of network training \cite{jagtap_adaptive_2020}. Adaptive weights leverage gradient statistics to optimize the interplay between various components in the loss function by incorporating additional weighting during the training process. While these solutions alleviate the convergence difficulties of neural networks in complex problems from different perspectives, there is insufficient evidence to demonstrate that simultaneously using different optimization schemes can still reduce network training difficulty and achieve large interval solutions. To solve the IVP problem of dynamic systems, Wen et al. proposed an ANN algorithm based on the Lie symmetry of differential equations \cite{Wen1, wen2, wen3, wen4}. This novel approach combines Lie group theory and neural network methods.

Given the aforementioned challenges, this paper proposes an innovative approach for obtaining large interval approximate solutions to IVPs of ODEs. In this algorithm, the interval is divided into several small compartments, and a neural network solution is learned on each compartment using PINN. Consequently, the trained ANNs generate fragments of the large interval solution of the IVP of an ODE on those sub-intervals, respectively. For each specific sub-interval, the ANN training requires neither complex structure nor dependence on large-scale training data, significantly reducing computational overhead. Subsequently, by assembling these neural network solutions, a piecewise expression of the large interval solution to the problem is constructed, termed the piecewise neural network solution. The compatibility of these sub-interval network solutions and the continuity and differentiability of the piecewise solution over the entire interval are investigated theoretically. Transfer learning techniques of network parameters between adjacent neural networks and multiple rounds of pre-training approaches in the training procedure are utilized.

The innovative contributions of this paper are threefold. First, a new ANN method for solving the extension problem of the local solution of the initial value problem of differential equations is presented. Second, a novel method of generalizing ANN and its application is introduced. Third, a new approach is explored in which modern algorithms (artificial neural networks, deep learning) are used to overcome the shortcomings of traditional methods for solving numerical solutions to differential equation problems. Our work aims to provide a more comprehensive understanding of the ANN method to solve ODE problems in both theoretical and practical aspects.

The structure of this paper is organized as follows. For completeness, Section \ref{sec:s2} briefly reviews the PINN method for IVPs of ODE systems; Section \ref{sec:s3} details the PWNN algorithm proposed in this paper; Section \ref{sec:s4} presents the theoretical analysis of PWNN solutions and implementation of the proposed method; Section \ref{sec:s5} demonstrates the applications of the proposed method in solving several specific IVPs of ODEs. Concurrently, comparisons between the results of our method and those of the PINN method and Runge-Kutta method are presented to demonstrate the efficiency of the proposed algorithm. Finally, Section \ref{sec:s6} summarizes the work.

\section{A brief recall of PINN}
\label{sec:s2}
Consider an IVP of a system of differential equations
\begin{equation}
    \label{eq:odes}
    \begin{cases}
        \dfrac{\mathrm{d}y}{\mathrm{d}x}=f(x, y), x \in I=(0, T), \\
        y(0)=y^0,
    \end{cases}
\end{equation}
with independent variable $x\in R$ and dependent variable $y=(y_1, y_2, \cdots, y_n)\in R^n$ and initial value (or initial condition, IC) $y^0=(y^0_1, y^0_2, \cdots, y^0_n)\in R^n$. We assume through the article that the right hand side function $f(x, y)=(f_1(x, y), f_2(x, y), \cdots, f_n(x, y)$ satisfies the condition \cite{ODE}

{\bf ($\textrm{H}_1$):} Continuous in $(x, y)\in \Omega\subset R\times R^n$ for a domain  $\Omega$ and Lipschitz continuous on $y$.

The condition guarantees the existence of unique local solution to the problem (\ref{eq:odes}). 

In addition, since the present article works on the large interval solution of IVP (\ref{eq:odes}), we also use the default assumption below:

{\bf ($\textrm{H}_2$):} The solution of IVP (\ref{eq:odes}) exists on large interval $I$ of $x$ (called large interval solution).

Thus, the PINN is an ANN algorithm for solving the solution of problem (\ref{eq:odes}) with fully connected structure that consists of one-dimensional input layer (0the layer), $M$ hidden layers, $s_i$ nonlinear active neurons for $i$th layer with $i=0, 1, \cdots, M+1$ and $s_0=1, s_{M+1}=n$, and output layer ($(M+1)$th layer) with $n$ linear output neurons. The output of the network is denoted as
\begin{equation*}
N(x)=(N_1(x), N_2(x), \cdots, N_n(x)), \label{netout}
\end{equation*}
for an input sample $x\in I$ which has mathematically the form of
\begin{equation}
N(x)=\sigma_M(\sigma_{M-1}(\cdots(\sigma_1(W^1 t+b^1)W^2+b^2)\cdots)W^M+b^M)W^{M+1}+b^{M+1},\label{trial}
\end{equation}
for selected activation functions $\sigma_{i+1}$, weights matrix $W^{i+1}$  between $i$th and $(i+1)$th layer and bias vector $b^{i+1}$ for $(i+1)$th layer neurons with $i=0, 1, 2, \ldots, M$ and $b_{M+1}(x)=x$.

The collection of all trainable parameters of the ANN, denoted as  $\vartheta$, is consisted of weights $W^i=(w^{i, 1}, w^{i, 2}, \ldots, w^{i, s_i})^T$ with $w^{i,j}=(w_i^{1,j}, w_i,^{2,j}, \cdots, w_i^{s_{i-1},j})$ in which the components $w^{k,j}_i$ ($k=1,2, \cdots, s_{i-1}, j=1, 2, \cdots, s_i$) are the weights connected the neurons in $(i-1)$th layer to $j$th neuron in $i$th layer, and layer biases $b^i=(b^{1}_i, b^{2}_i, \ldots, b^{s_i}_i)^T (i=1, 2, \cdots, M+1)$.  Hence, $W^i$ is a $s_{i}\times s_{i-1}$ matrix. Let $\vartheta$ denote the set of all weights and biases parameters of the net work. The unsupervised learning (training) method is used. That is, the optimal solution $N(x)$ at optimal parameters $\vartheta=\vartheta^*$ is searched such that the loss function to be given below is minimized.

\begin{figure}[htbp]
    \centering
    \includegraphics[width=0.9\textwidth]{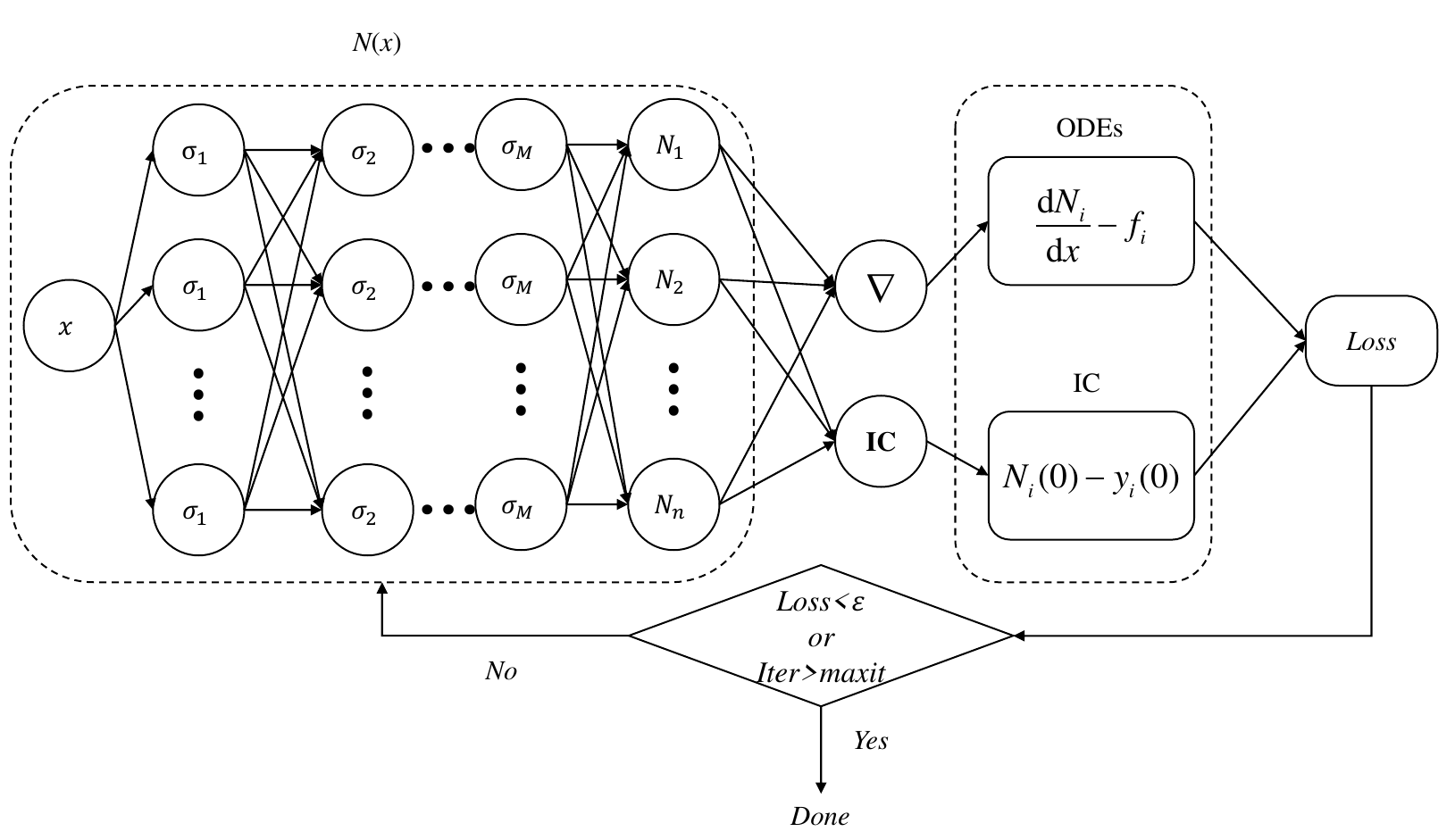}
    \caption{\centering Network structure of PINN for solving IVP (\ref{eq:odes}).}
    \label{fig:pinn_odes}
\end{figure}

Usually, the mean square error (MSE) is used to define the loss function for training the network as follows:
\begin{equation}
    \label{eq:loss}
    Loss = \frac{1}{n}\sum_{i=1}^n MSE_{f_i}+\frac{1}{n}\sum_{i=1}^n (N_i(0)-y_i^0)^2,
\end{equation}
where
\begin{equation*}
    MSE_{f_i} = \frac{1}{M}\sum_{j=1}^M (N_i'(x_j)-f_i(x_j, N_1(x_j), N_2(x_j),\cdots, N_n(x_j)))^2,
\end{equation*}
and $\{x_j\}_{j=1}^M$ is a training data set sampled by a distribution sense from the solution interval $I$, and $N_i'(x)$ is the derivative of the approximate solution $N_i(x)$ with respect to the independent variable $x$. In training the network, the automatic differentiation operation $\nabla$ is used to realize $N_i'(x)$. One of initialization methods \cite{glorot_xavier_2010} is performed on the neural network parameters when starting the training. A optimization technique \cite{kingma_adam_2017} is used to adjust the network parameters in order to minimize the loss function (\ref{eq:loss}) and carry out the back propagation of errors.
After completing the training of the network by minimizing the loss function ($Loss<\epsilon$) or sufficient iterative parameter refinement (learning, training) ($Iter>maxit$), we obtain the network pattern (\ref{netout}) with a optimal parameter set $\vartheta^*$, which represents the approximate solution of problem (\ref{eq:odes}). The basic structure of above mentioned ANN (e.g, PINN) is shown in Fig. (\ref{fig:pinn_odes}). 

Particularly, what should be especially mentioned here that the experiences of applying an ANN on solving IVP of a (partial) differential equations show that it often yields high accuracy approximate solution merely near the initial value point rather than the entire interval $[0, T)$ \cite{ODE}.
That is, generally, there exist a constant $0<\delta\ll T$ such that the output $N(x)$, called PINN solution, converges to exact solution $y(x)$ over $[0, \delta)$ with high accuracy. We denote the approximation as
\begin{equation*}
N(x)\approx y(x).\label{appr}
\end{equation*}
In other words, as the training point gradually moves away from the starting point $x=0$, the learning ability or generalization ability of the ANN gradually weakens. There are two main potential reasons for this phenomenon. First, the requirement of $Loss\to 0$ in network training is only a necessary condition for $N(x)$ to solve IVP (\ref{eq:odes}). Second, with the expansion of the sampling data interval, the support of initial value information to the training ANN may gradually decline \cite{lagaris_artificial_1998}.
Therefore, in order to obtain the large interval solution of problem (\ref{eq:odes}) by ANN method, it is often necessary to construct a neural network with a large number of neurons and a deep layers so that the influence of the initial value is deep and wide, which may make the network learn the solution as accuracy as possible \cite{NUrODE}. However, the increase of network scale significantly increases the training difficulty of the neural network, and even makes it impossible to get normal training.
In view of this problem, based on the advantage of PINN's high-precision learning ability in a neighborhood of initial value point, we propose a PWNN method in next section to obtain the large interval solution. In the method, the solution interval $I$ is split into several intercell compartments and a related small size in scale neural network is applied on each the sub-interval for solving a sub-IVP corresponding to (\ref{eq:odes}). Finally stitches these neural network solutions of the sub-intervals together to obtain the approximate solution of problem (\ref{eq:odes}) over whole interval $I$.

\section{A piecewise ANN method}
\label{sec:s3}
In this section, based on the default assumptions ({\rm $\textrm{H}_1)-(\textrm{H}_2$}) and well learning ability of a PINN, we will present the PWNN method for finding the large interval approximate solution of IVP (\ref{eq:odes}).
\subsection{Method}

The basic framework of the method is described as follows.

\begin{figure}[htbp]
    \centering
    \includegraphics[width=0.6\textwidth]{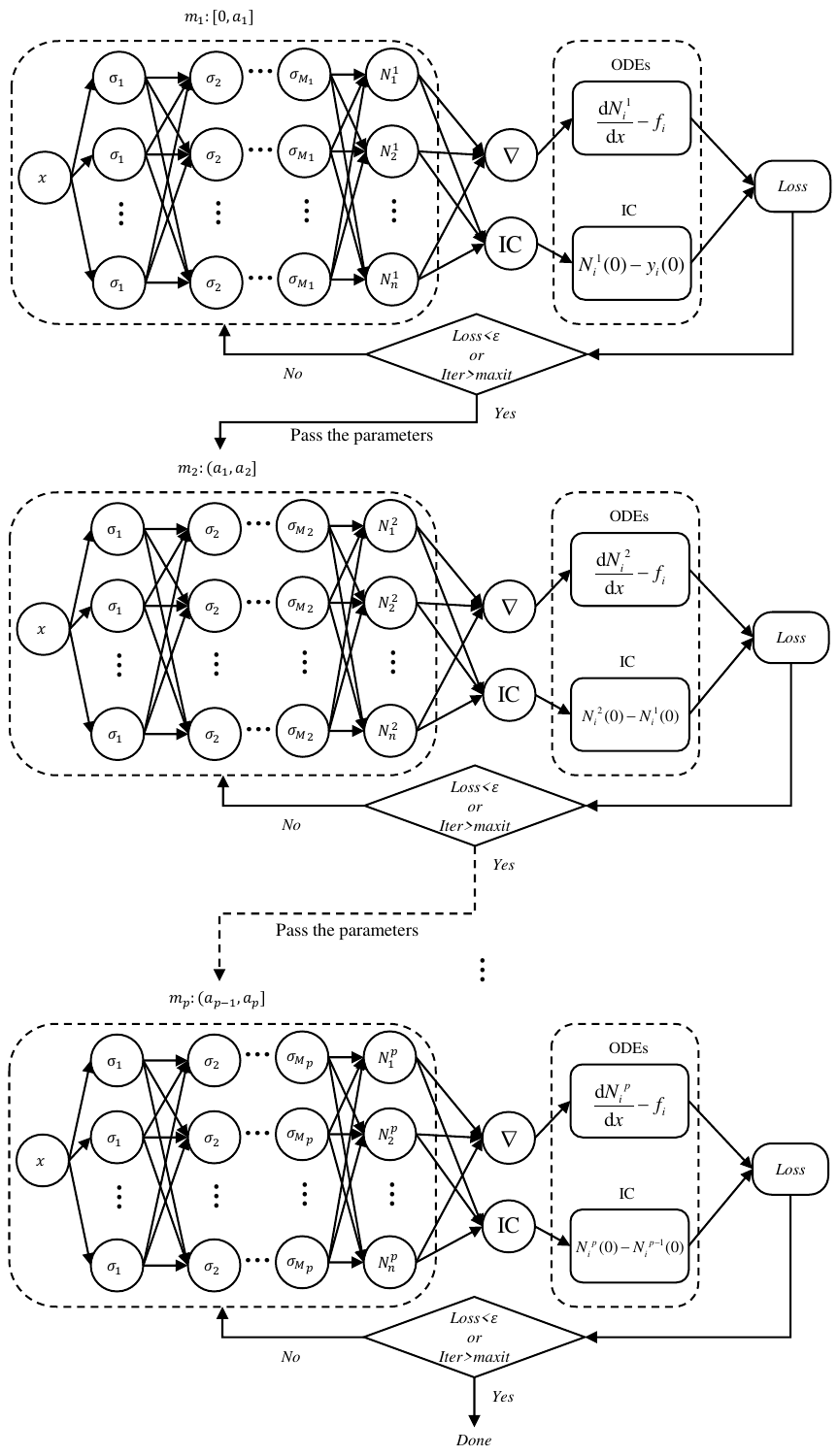}
    \caption{\centering Network structure of PWNN for solving IVP (\ref{eq:odes}).}
   \label{fig:fdwlyuanli}
\end{figure}

First, we divide the interval $[0, T]$ into $p$ sub-intervals, namely insert $p-1$ points in interval $[0, T]$ making
\[0 = a_0 < a_1 < \cdots < a_p = T, \,\, [0, T] = \bigcup_{i=1}^p \Delta_i, \,\, \Delta_i=[a_{i-1}, a_i],\]
and suppose that sub-interval $\Delta_i$ is relative small so that the well performing of PINN is guaranteed.
Then, we construct $p$ ANNs (e.g. PINNs) with identical structures as given in Section 2 corresponding to the splitting of the interval $[0, T]$, denoted them as $\{m_i\}_{i=1}^p$.
Each neural network $m_i$ is used to solve a sub-IVP of the same differential equations in (\ref{eq:odes}) on the $i$th sub-interval $\Delta_i$ with a well defined initial value at $a_{i-1}$.
The output of the network is denoted as $N^{i}(x)=(N_1^{i}(x), N_2^{i}(x), \cdots, N_n^{i}(x))$, which has composite expression like (\ref{trial}) and represents the approximation of exact solution of the sub-IVP, denoted as $\bar{y}^i(x)=(\bar{y}^i_1(x), \bar{y}^i_2(x), \cdots, \bar{y}_n^i(x))$. Particularly, component $N_j^{i}(x)$ is the output of $j$-th neuron of network $m_i$ in the output layer for an input $x \in \Delta_i$, representing the approximation of $\bar{y}^i_j(x)$ for $j \in \{1, 2, \cdots, n\}$. The network is illustrated in Fig. \ref{fig:fdwlyuanli}.

Specifically, first, a PINN $m_1$ is trained to solve the sub-IVP of the system in (\ref{eq:odes}) on interval $\Delta_1$ with initial value $\bar{y}^1(0)=y^0$. The loss function used for the network is defined as
\begin{align*}
    Loss_1 &= \frac{1}{M_1 n}\sum_{i=1}^n \sum_{j=1}^{M_1} (\dfrac{\mathrm{d}N_i^{1}}{\mathrm{d}x}(x_j)-f_i(x_j, {N_1}^{1}(x_j),{N_2}^{1}(x_j), \cdots, N_n^{1}(x_j)))^2\\
            &+ \frac{1}{n}\sum_{i=1}^n (N_i^{1}(0)-y_{i}^0)^2.
\end{align*}
Here, the training data points $\{x_j\}_{j=1}^{M_1}$ are taken from the interval $\Delta_1$ by some distribution sense (the following-up data samplings are similar). The training of the network $m_1$ is completed when the termination (convergence) conditions of the network training (generally the loss function value drops to near zero or the iteration step size is sufficient) are satisfied. As a result, due to the approximation property of ANN, the ANN is an approximate solution of the sub-IVP, denoted by $N^{1}(x)=(N_1^{1}(x), N_2^{1}(x), \cdots, N_n^{1}(x))$, i.e., it is the approximation of exact solution $\bar{y}^1(x)=(\bar{y}^1_1(x), \bar{y}^1_2(x), \cdots, \bar{y}^1_n(x))$ of the sub-IVP over $\Delta_1$. Notably, compared to a neural network that solves IVP (1) over the entire interval $I$, network $m_1$ requires only a structure of relatively small scale and is trained well on $\Delta_1$ because it is a small neighborhood of the initial value point $x = 0$. This gives full play to the strong local learning ability, faster convergence speed and high approximation accuracy of PINN.

After obtaining the first network $m_1$ on $\Delta_1$, we train and use the next network $m_2$ to approximately solve the sub-IVP of the system in (\ref{eq:odes}) with an initial value of $\bar{y}^2(a_1)=N^1(a_1)$ in terms of previous trained network $m_1$ on the interval $\Delta_2$. The corresponding loss function for $m_2$ is defined as follows:
\begin{equation}
    \begin{aligned}
        \label{eq:m2loss}
        Loss_2
        &= \frac{1}{M_2 n}\sum_{i=1}^n \sum_{j=1}^{M_2} \left(\frac{\mathrm{d}N_i^{2}}{\mathrm{d}x}(x_j)-f_i(x_j,
          {N_1}^{2}(x_j), {N_2}^{2}(x_j), \cdots, N_n^{2}(x_j))\right)^2\\
        &+ \frac{1}{n}\sum_{i=1}^n \left(N_i^{2}(a_1)-
          N_i^{1}(a_1)\right)^2.
    \end{aligned}
\end{equation}
Here, the training data points $\{x_j\}_{j=1}^{M_2}$ are taken from the interval $\Delta_2$. The second term on the right-hand side of the loss function ensures that network $m_2$ approximates $m_1$ as closely as possible at the initial point $a_1$. When the termination conditions are satisfied network $m_2$ has been trained and yields approximate solution $N^{2}(x)=(N_1^{2}(x), N_2^{2}(x), \cdots, N_n^{2}(x))$ of the exact solution $\bar{y}^2(x)=(\bar{y}^2_1(x), \bar{y}^2_2(x), \cdots, \bar{y}^2_n(x))$ of the sub-IVP.

Inductively, suppose we have already trained network $m_{i}$ and obtained approximate solution $N^{i}(x)=(N_1^{i}(x), N_2^{i}(x), \cdots, N_n^{i}(x))$ of the exact solution $\bar{y}^i(x)=(\bar{y}^i_1(x), \bar{y}^i_2(x),$ $ \cdots, \bar{y}^i_n(x))$ of the sub-IVP of system in (\ref{eq:odes}) with initial value $\bar{y}^{i}(a_{i})=N^{i-1}(a_{i-1})$ on the interval $\Delta_i$ for $i=1, 2, \cdots, k-1$. Now, next network $m_k$ is trained and used to approximate exact solution $\bar{y}^k (x)$ of the sub-IVP of the same system in (\ref{eq:odes}) on the interval $\Delta_k$ with initial value $\bar{y}^k(a_{k-1})=N^{k-1}(a_{k-1})$ in terms of previous network $m_{k-1}$. The training data points $\{x_j\}_{j=1}^{M_k}$ are taken from the interval $\Delta_k$. The corresponding loss function for network $m_k$ is defined as:
\begin{equation}
\label{eq:mkloss}
\begin{aligned}
    Loss_k &= \frac{1}{M_k n}\sum_{i=1}^n \sum_{j=1}^{M_k} \left(\frac{\mathrm{d}N_i^{k}}{\mathrm{d}x}(x_j)-f_i(x_j, {N_1}^{k}(x_j), {N_2}^{k}(x_j), \cdots, N_n^{k}(x_j))\right)^2 \\
    &+ \frac{1}{n}\sum_{i=1}^n \left(N_i^{k}(a_{k-1})-N_i^{{k-1}}(a_{k-1})\right)^2.
\end{aligned}
\end{equation}
When the termination conditions are satisfied network $m_k$ has been trained and yields approximate solution $N^{k}(x)=(N_1^{k}(x), N_2^{k}(x), \cdots, N_n^{k}(x))$ of the exact solution $\bar{y}^k(x)=(\bar{y}^k_1(x), \bar{y}^k_2(x), \cdots, \bar{y}^k_n(x))$ of the sub-IVP.

In this way we train and use all ANNs $m_i$ or equivalently $N^i(x)$, subsequently to approximate exact solution $\bar{y}^i(x)$ of the sub-IVP over sub-intervals $\Delta_i$ with initial value $\bar{y}^i(a_{i-1})=N^{i-1}(a_{i-1})$ for $i=1, 2, \cdots, p$.

Finally, we construct a piecewise function over $[0, T)$ by these PINN's solutions given by:
\begin{equation*}
    \hat{y}(x) = \begin{cases}
        N^1(x), & x \in [0, \,\, a_1], \\
        N^2(x), & x \in (a_1, \, \, a_2], \\
        \vdots \\
        N^p(x), & x \in (a_{p-1}, T),
    \end{cases}
\end{equation*}
where
\[N^k(x) = (N_1^{k}(x), N_2^{k}(x), \cdots, N_n^{k}(x)), \,\, k=1, 2, \cdots, p.\]

\section{Theoretical analysis and a parameter transfer method}
\label{sec:s4}
In the section, we will theoretically demonstrate that the function $\hat{y}(x)$ given in last section is an approximation of the large interval solution $y(x)$ of IVP (\ref{eq:odes}) we are looking for. Meantime, we give a parameter transfer method and multiple rounds of pre-training approach in training PWNN.
\subsection{Approximation of large interval solution}

In fact, the network $m_k$, or equivalently $N^k(x)$ with $N^0(x)=y^0$, presented in the last section approximately solved the following sub-IVP of the system in (\ref{eq:odes}):
\begin{equation}
\begin{cases}
\dfrac{\mathrm{d}\bar{y}^k}{\mathrm{d}x} = f(x, \bar{y}^k), \,\,x\in \Delta_k, \\
\bar{y}^k(a_{k-1}) = N^{k-1}(a_{k-1}),
\end{cases}\label{ybar}
\end{equation}
for dependent variable $\bar{y}^k(x)=(\bar{y}^k_1(x), \bar{y}^k_2(x), \cdots, \bar{y}^k_n(x))$ and for each $k\in\{1, 2, \ldots, p\}$, inductively, and we have
\begin{equation}
N^k(x)\approx \bar{y}^k(x) \,\, \textrm{ on } \,\, \Delta_k, \,\, k=1, 2, \cdots, p. \label{Nbary}
\end{equation}

In addition, we also set a sub-IVP of (\ref{eq:odes}) on each sub-interval $\Delta_k$ as
\begin{equation}
    \begin{cases}
        \dfrac{\mathrm{d}y^k}{\mathrm{d}x}=f(x, y^k), x \in \Delta_k, \\
        y^k(a_{k-1})=y^{k-1}(a_{k-1}),
    \end{cases} \label{nobary}
\end{equation}
with $y^k(x)=(y^k_1(x), y^k_2(x), \cdots, y^k_n(x))$ and $y^0(0)=y^0$ for $k\in \{1, 2, \cdots, p\}$. Due to the uniqueness of the solution of the IVP of an ODEs under assumption ($\textrm{H}_1$) given in last section, we know that solution $y^k(x)$ of sub-IVP (\ref{nobary}) is the restriction of the exact solution $y(x)$ of IVP (\ref{eq:odes}) on $\Delta_k$. That is,
\begin{equation}
y^k(x)=y(x)\,\, \textrm{ on } \,\,  \Delta_k, \,\, k=1, 2, \cdots, p, \label{yeqyk}
\end{equation}
and particularly,
$
y(x)=y^1(x)=\bar{y}^1(x) \,\, \textrm{ on } \,\, \Delta_1.
$
Hence, noticing (\ref{Nbary}), approximation
\begin{equation}
  N^1(x)\approx y(x), \label{N1y1}
\end{equation}
holds on $\Delta_1$.

For $k=2$, we have
$$y^2(a_1)\overset{(\ref{nobary})}{=}y^1(a_1)=\bar{y}^1(a_1) \overset{(\ref{Nbary})}{\approx}N^1(a_1)\overset{(\ref{ybar})}{=}\bar{y}^2(a_1).$$
Hence, on sub-interval $\Delta_2$ the initial values of IVP (\ref{ybar}) and IVP (\ref{nobary}) are approximate. Consequently, by the qualitative property that the solution of IVP of an ODEs with assumption ($\textrm{H}_1$) continuously depends on the initial value, approximation
\begin{equation}
N^2(x)\approx y(x), \label{N2y2}
\end{equation}
holds on $\Delta_2$.

Similarly, we can inductively prove
\begin{equation}N^k(x)\approx y(x), \label{Nkyk}\end{equation}
hold on intervals $\Delta_k$ for $k=1, 2, \cdots, p$. This shows that approximation $\hat{y}(x)\approx y(x)$ holds on each $\Delta_k$.

Furthermore, from the construction of PINN $N(x)$ (see the multiple composite structure of the PINN's output expression (\ref{trial})), we see that function $N^k(x)$ is continuous and differential over the each sub-interval $\Delta_k$. Therefore, except the end points $a_k$ of $\Delta_k (k=1, 2, \cdots, p-1)$, function $\hat{y}(x)$ is continuous and differential over large interval $[0, T)$. Although we do not confirm the continuity of $\hat{y}(x)$ on thoses end points, we, from (\ref{ybar}) and (\ref{Nbary}), know approximation
$$N^k(a_k)\approx N^{k+1}(a_k),$$
implying the approximation of the two side limitations
\begin{equation}
\hat{y}(a_k^-)\approx \hat{y}(a_k^+),\label{endpoint}
\end{equation}
hold for $k=1, 2, \cdots, p-1$. This indicates that there may exist just a small 'jump' in the value of function $\hat{y}(x)$ at an end point. However, in numerical solution, the errors can be controlled within the tolerable error range by improving the convergence accuracy of PINN solutions $N^k(x)$.

To sum up above procedure, we have proved the following conclusions.

{\bf Theorem: }On the solution interval $[0,T)$, the piecewise function $\hat{y}(x)$ satisfies the approximation $\hat{y}(x)\approx y(x)$ and is continuously differentiable except for the finite points $a_k(k=1, 2, \cdots, p-1)$, where the approximation (\ref{endpoint}) holds.

The Theorem is the theoretical basis of our PWNN method.
\subsection{Transfer of network parameters}

According to the research by Xavier et al., the parameter initialization of a neural network directly affect whether the network can be trained successfully \cite{glorot_xavier_2010}. In the proposed PWNN method, we need to train $p$ interrelated PINNs. Therefore, the parameter initialization and transfer of the PWNNs in training them are a key step to successfully obtain large interval approximate solution $\hat{y}(x)$. In order to assess and improve the stability and performance in certain cases, there are some additional network parameter transfer and multiple rounds of pre-training techniques which we employ beyond the basic setup.

1. A parameter transfer technique. 
In fact, we have already used a parameter initialization technique in the training and theoretical analysis of PWNN. This is shown in IVP (\ref{ybar}) and procedure of constructing ANNs $m_k (k=1, 2, \cdots, p)$, which is explained in more details in the following steps.

Start with $k=1$. By the standard steps of PINN, initializing network parameters $\vartheta_1$ with Xavier technique or some other ones, we then train first PINN $m_1$ on $\Delta_1$ to solve IVP (\ref{ybar}) for $k=1$ and determine $N^1(x)$.

For $k=2, 3, \cdots, p$, initializing parameters $\vartheta_k$ of $m_k$ using parameters $\vartheta_{k-1}$ obtained in $m_{k-1}$, we then train the $k$th PINN $m_k$ on $\Delta_{k}$ to solve IVP (\ref{ybar}) and determine $N^k(x)$. That is, once network $m_{k-1}$ has completed training, we pass parameters $\vartheta_{k-1}$ to network $m_k$ as initialization of parameters $\vartheta_k$ of $m_k$ as shown in Fig.2.

It is worth noting here that since network $m_{k}$ inherits the parameters of $m_{k-1}$, the computational effort of loss function $Loss_k$ will be greatly reduced, thus speeding up the training of the network.

2. Multiple rounds of pre-training method. In the first time training PINN, it is not always satisfactory performance of the being trained network due to the randomly parameter initialization. In this case, it is common to conduct training the network in several rounds (pre-training) based on the latest obtained parameters. That is, the obtained parameters after completing the training of an ANN are used as initialization of the network to train again, and so on until the end of the training round. This multi-round pre-training gradually optimizes the parameter initialization, avoids the uncertainty caused by random initialization and guides the improving training quality. 

In our PWNN case, the above parameter transfer method (see above the case $k=1$) is used in the first round of training, and then the multi-round pre-training method is used after the second round. Specifically, each neural network $m_j, j= 2, 3, \cdots, p$ receives the parameters of $m_{j-1}$ as parameter initialization in the first round of training. For subsequent rounds, in training a network, it receives the network parameters obtained from the latest training round of this network as its initialization. This multi-round pre-training will progressively improve the approximation of $\hat{y}(x)$ to the exact solution of IVP (\ref{eq:odes}). To further
improve the approximation accuracy of $\hat{y}(x)$ to exact solution $y(x)$, we even apply the two kind techniques interactively in training PWNNs.
\subsection{Implementation of PWNN}
Let $\vartheta_k^i$ denote the training parameter set of network $m_k$ in the $i$th round and $maxit$ denote the maximum number of iterations. $Loss$ represents the loss value of the neural network, which is determined by Eq.(\ref{eq:mkloss}), and $\varepsilon>0$ represents a pre-specified error limitation.

Now, summarizing the above statement of PWNN method, we have the following {\bf Algorithm PWNN} as an implementation of our proposed PWNN method.

\begin{algorithm} 
   \renewcommand{\thealgorithm}{}
	\caption{PWNN} 
	\begin{algorithmic}
        \State {\bf Initialization:} Define a generic functional module
        \State Back-propagation=
        \State\hspace{0.5cm}\{
        \State \hspace{1cm}{\bf While} {$it\leq maxit \land loss \geq \varepsilon$}
        \State \hspace{1.5cm}Update $\vartheta_k^i$ via a optimizer \cite{kingma_adam_2017}
        \State \hspace{1.5cm}Update the loss based on Eq.(\ref{eq:mkloss}) and back propagation method
        
        \State\hspace{0.5cm}\}
        \State {\bf Begin:}
        \For{$i=1$ to $n$}

        \If{$i=1$}
        \For{$k=1$ to p}

        \If{$k=1$}
        \State $\vartheta^i_k$ is initialized using Xavier method or some others \cite{glorot_xavier_2010}
        \Else 
        \State initializing $\vartheta_k^i=\vartheta_{k-1}^{i}$
        \EndIf

        \State \hspace{0.5cm}Back-propagation
        \EndFor

        \Else
        \State initializing $\vartheta_k^i=\vartheta_k^{i-1}$
        
        \State Back-propagation
        
        \EndIf
        
        \EndFor
        \State {\bf End}
	\end{algorithmic} 
\end{algorithm} 

\section{Experiment}
\label{sec:s5}
In this section, we give several numerical experiments using the PWNN algorithm, and compare the results with those of  PINN method and Runge-Kutta method to show the validity of our proposed method.

\subsection{Example 1}
\begin{figure}[ht]
    \centering
    \includegraphics[width=0.9\textwidth]{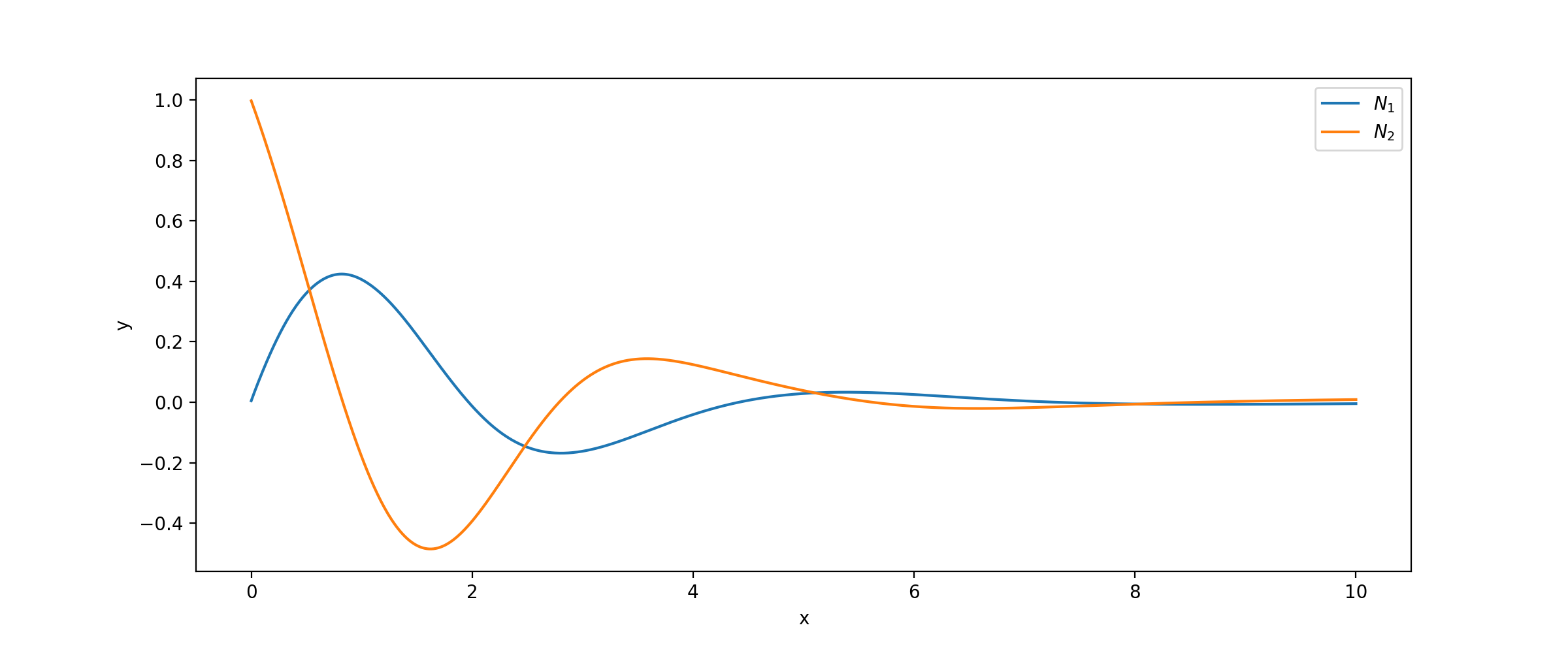}
    \caption{\centering PINN results for IVP (\ref{eq:exam1}) with final loss function value $6.58*10^{-5}$.}
    \label{fig:NN_exam1}
\end{figure}

Consider the following IVP of an ODEs
\begin{equation}
    \label{eq:exam1}
    \begin{cases}
        \frac{\mathrm{d}y_1}{\mathrm{d}x}=y_2,\\
        \frac{\mathrm{d}y_2}{\mathrm{d}x}=-y_2-(2+\sin x)y_1,\\
        y_1(0)=0, y_2(0)=1,
    \end{cases}
\end{equation}
with dependent variable $y(x)=(y_1(x), y_2(x))$ on interval $[0,10]$.
We first use normal PINN to find the ANN solution of the problem. We construct a 3-layer neural network with 1 input neuron, 2 output neurons, and 2 hidden layers of 20 neurons each.
We uniformly sampled 1,000 points on the interval $[0, 10]$ as the training data set. The learning rate is set to $0.01$, and the number of iterations is 10,000. The training results of the network are shown in Fig. \ref{fig:NN_exam1}.

\begin{figure}[htbp]
    \centering
    \includegraphics[width=0.9\textwidth]{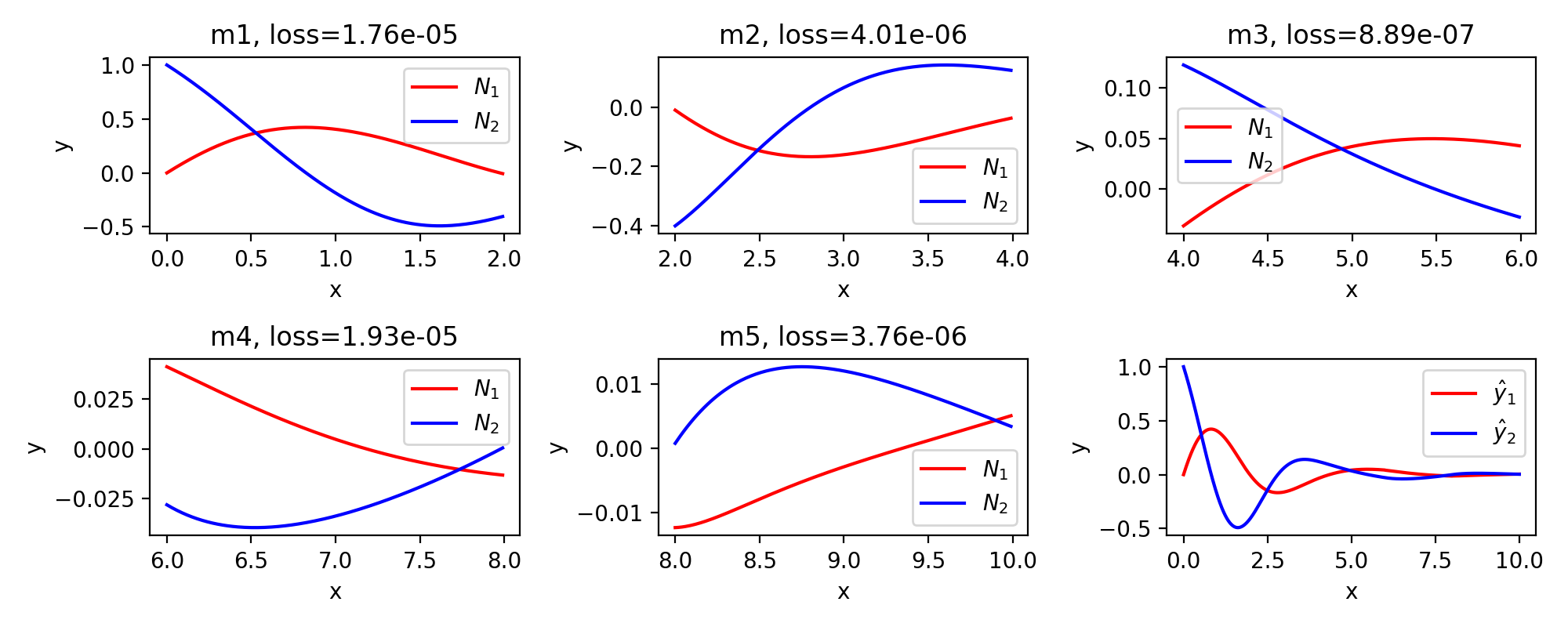}
    \caption{\centering Result figures of PWNNs: the first five is $m_k$ with final loss function value on each sub-interval $\Delta_k (k=1, 2, \cdots, 5)$. The last one is their combinations $\hat{y}_1(x)$ and $\hat{y}_2(x)$ .}
    \label{fig:fdw_exam1}
\end{figure}

Then we use PWNN to solve the problem with less training data and the same structure as PINN above on each sub-interval. Interval $[0, 10]$ is divided into five equidistance sub-intervals $\Delta_i=[2(i-1), 2i] (i=1, 2, \cdots, 5)$.
The neural network employed in each sub-interval consists of two hidden layers, with each hidden layer containing 20 neurons. In order to compare it with PINN, we uniformly sample 200 data points from each sub-interval as the training data. The learning rate is set to 0.01, and each neural network performed 2000 iterations. The results of training PWNNs $m_i (i=1, 2, \cdots, 5)$ are shown in Fig. \ref{fig:fdw_exam1}. The first five graphs represents the PWNNs on each interval as well as loss function values, and the last graph is their combination $\hat{y}(x)=(\hat{y}_1(x), \hat{y}_2(x))$.

We also solve the problem using the fourth-order Runge-Kutta method to compare the results with those of using above methods.

\begin{figure}[htbp]
    \centering
    \includegraphics[width=0.9\textwidth]{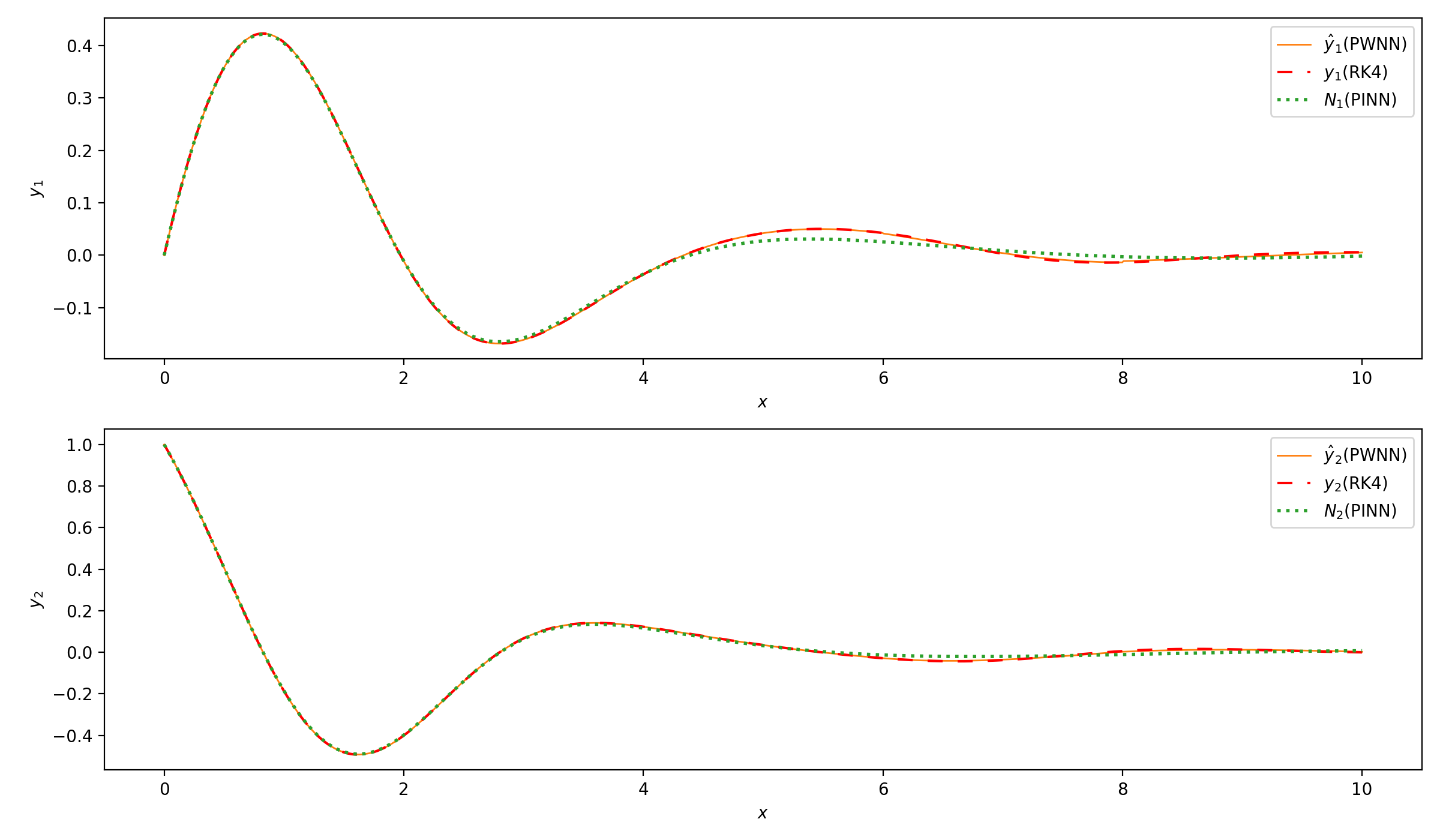}
    \caption{\centering Comparisons of the results of solving IVP (\ref{eq:exam1}) using PINN, PWNN and RK4 methods.}
    \label{fig:Three_graph_comparison}
\end{figure}

In Fig. \ref{fig:Three_graph_comparison}, the results obtained by using PINN method, PWNN approach and fourth-order Runge-Kutta method (RK4) are  presented. It can be observed that the overall results of the three methods are similar. However, compared to PINN, the results of PWNN and RK4 are closer.
It shows that PWNN is more efficiency than PINN in the sense of numerical solution.

\subsection{Example 2}
\label{sec:4.2}
\begin{figure}[htbp]
    \centering
    \includegraphics[width=0.9\textwidth]{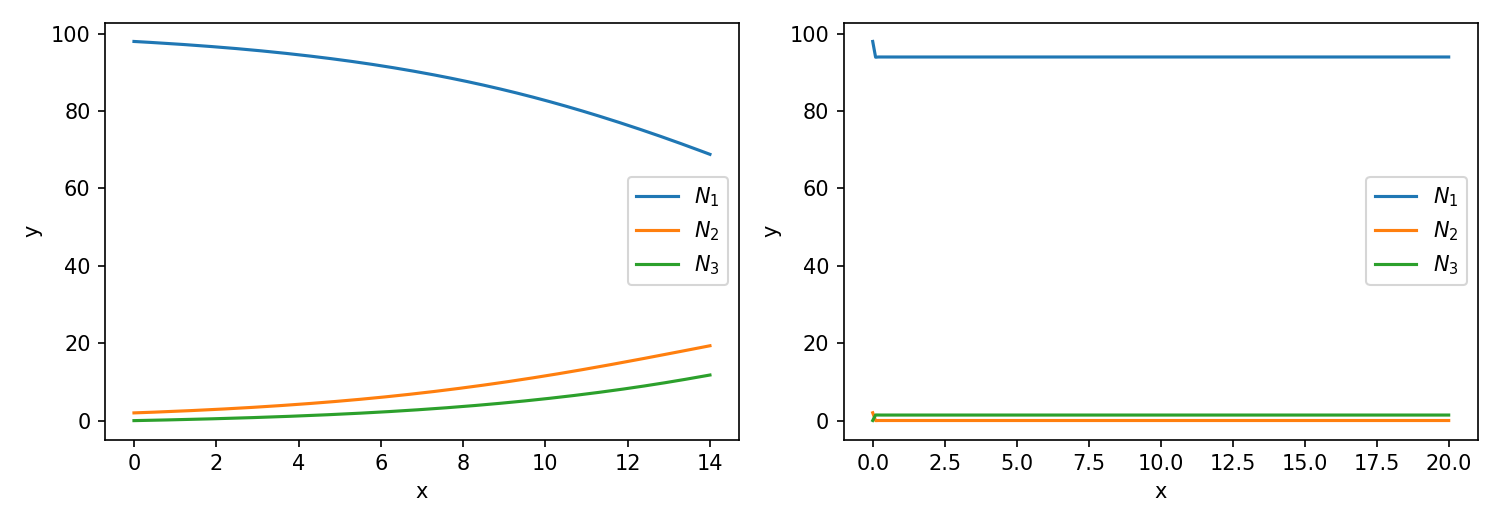}
    \caption{\centering Results of PINN to solve IVP (\ref{eq:sir}) which failed to get large interval solutions.}
    \label{fig:fig1}
\end{figure}

Consider the following IVP of a SIR model of epidemic dynamics
\begin{equation}
    \label{eq:sir}
    \begin{cases}
        \frac{\mathrm{d}y_1}{\mathrm{d}x} = -0.003 y_1 y_2,\\
        \frac{\mathrm{d}y_2}{\mathrm{d}x} = 0.003y_1 y_2 - 0.1 y_2,\\
        \frac{\mathrm{d}y_3}{\mathrm{d}x} = 0.1y_2,\\
        y_1(0) = 98,y_2(0)=2,y_3(0)=0.
    \end{cases}
\end{equation}
First, we use PINN to find its solution over the interval $[0, 50]$. The used PINN consists of $1$ input neuron, $3$ hidden layers of $20$ neurons each, and $3$ output neuron. The 200 training points are taken from $[0, 20]$ equidistantly.
The network adopts Xavier initialization, employs the Adam optimization algorithm with a learning rate $0.01$, and runs for 10000 iterations. Fig. \ref{fig:fig1} illustrates the solutions given by the PINN for IVP (\ref{eq:sir}) over the interval $[0, 20]$. Observably, the PINN successfully computes solutions merely within the intervals $[0, 14]$ 
 (left graph) instead of large interval $[0, 20]$ (right graph). That is, as the solution interval further extends the learning results of PINN notably deviate from the real case.

We now solve the problem over a large interval $[0,50]$ using the PWNN method with a smaller size structure than the PINN above. Divide interval $[0, 50]$ into $10$ equally sized sub-intervals and Take 100 training points in each sub-interval by uniform distribution. Construct a neural network for each segment, where each network consists of $1$ input neuron, $1$ hidden layer with $20$ neurons, and $3$ output neurons. Utilize Xavier initialization and Adam optimization algorithms for each network with learning rate $0.01$ and 10000 iterations to train $m_1$ on  interval $[0, 5]$.
\begin{figure}[htbp]
    \centering
    \includegraphics[width=0.9\textwidth]{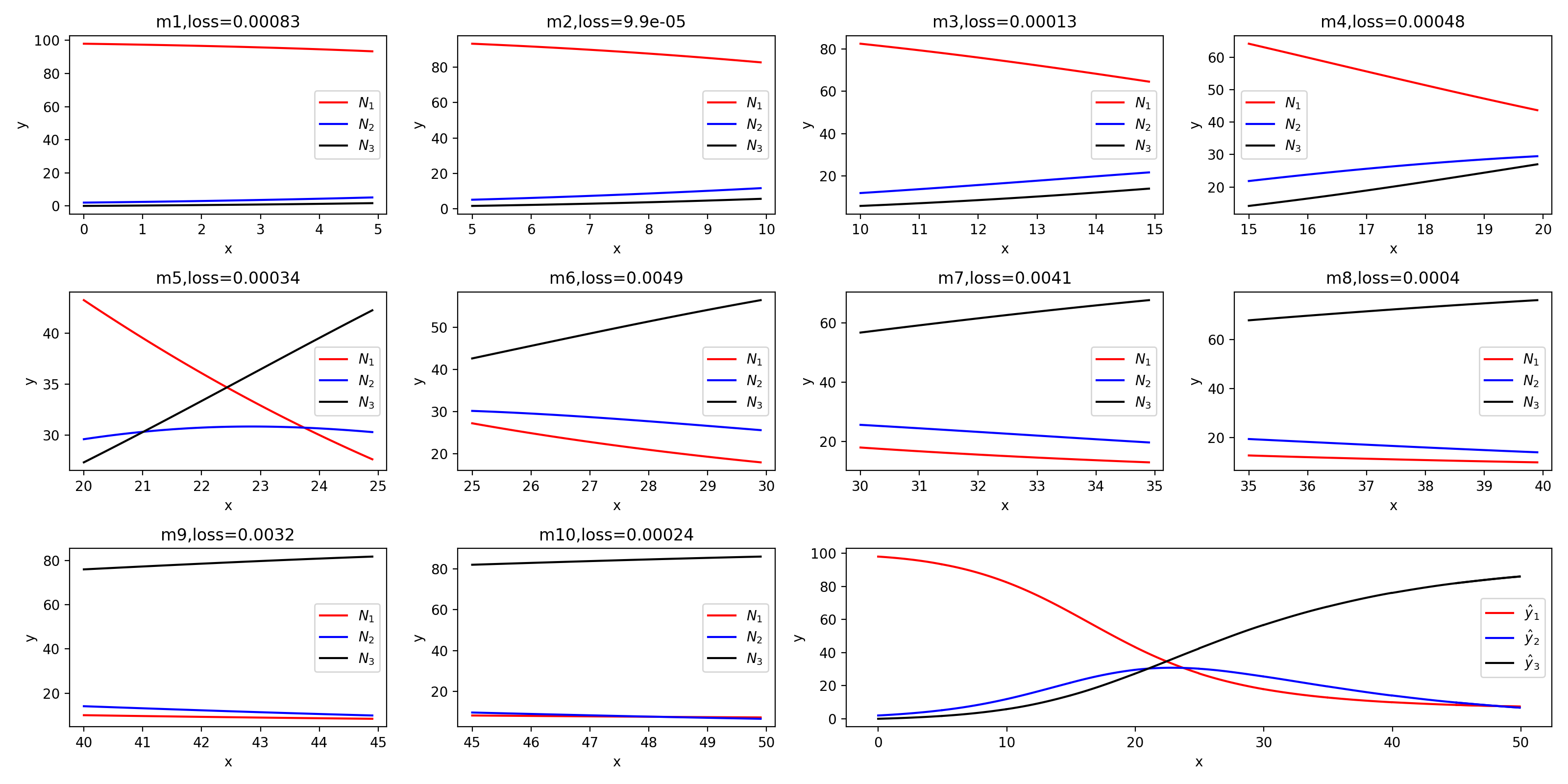}
    \caption{\centering Results of PWNNs $m_i$ on $\Delta_i$ for $(i=1, 2, \cdots, 10)$ and loss function values for each PWNN.}
    \label{fig:fig2}
\end{figure}
Fig. \ref{fig:fig2} illustrates the training results of PWNNs for all sub-intervals. The titles of the figures specify the names of the networks for each interval along with the respective values of the loss function after the completion of network training.
Fig. \ref{fig:fig3} displays the comparisons between the results of both PWNNs and the RK4 method. The solid lines represent the solutions $\hat{y}(x)$ obtained from PWNNs, while the discrete points depict the solutions obtained using RK4 method. Obviously, PWNN gives a large interval solution of problem (\ref{eq:sir}) that is highly consistent with RK4.
\begin{figure}[htbp]
    \centering
   \includegraphics[width=0.9\textwidth]{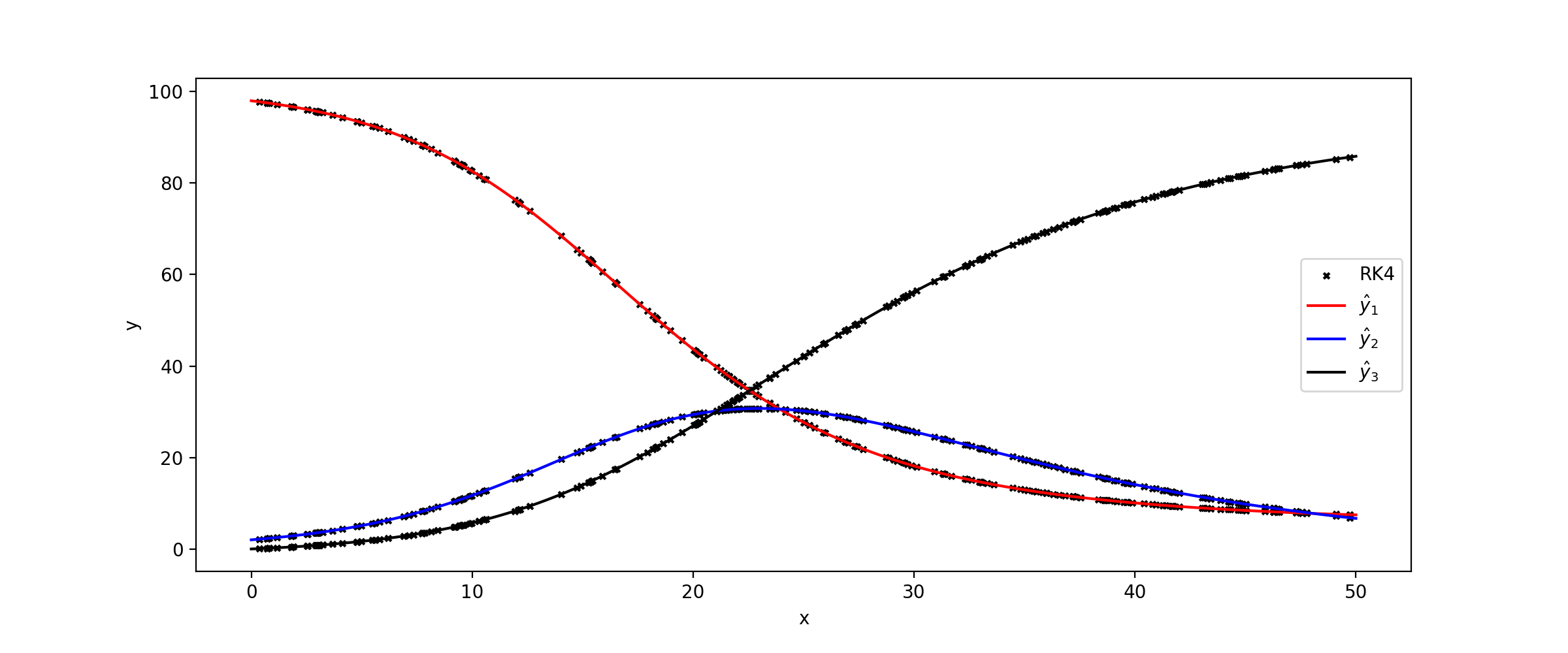}
   \caption{\centering Comparison of the results of using PWNN and RK4 methods.}
   \label{fig:fig3}
\end{figure}

Example 2 shows that direct application of conventional PINN can only obtain the inter-cell solution (local solution) of the problem, while PWNN directly gives a large interval solution consistent with RK4, which demonstrates the effectiveness of the given PWNN algorithm.

\subsection{Example 3}
\begin{figure}[htbp]
   \centering
   \includegraphics[width=0.9\textwidth]{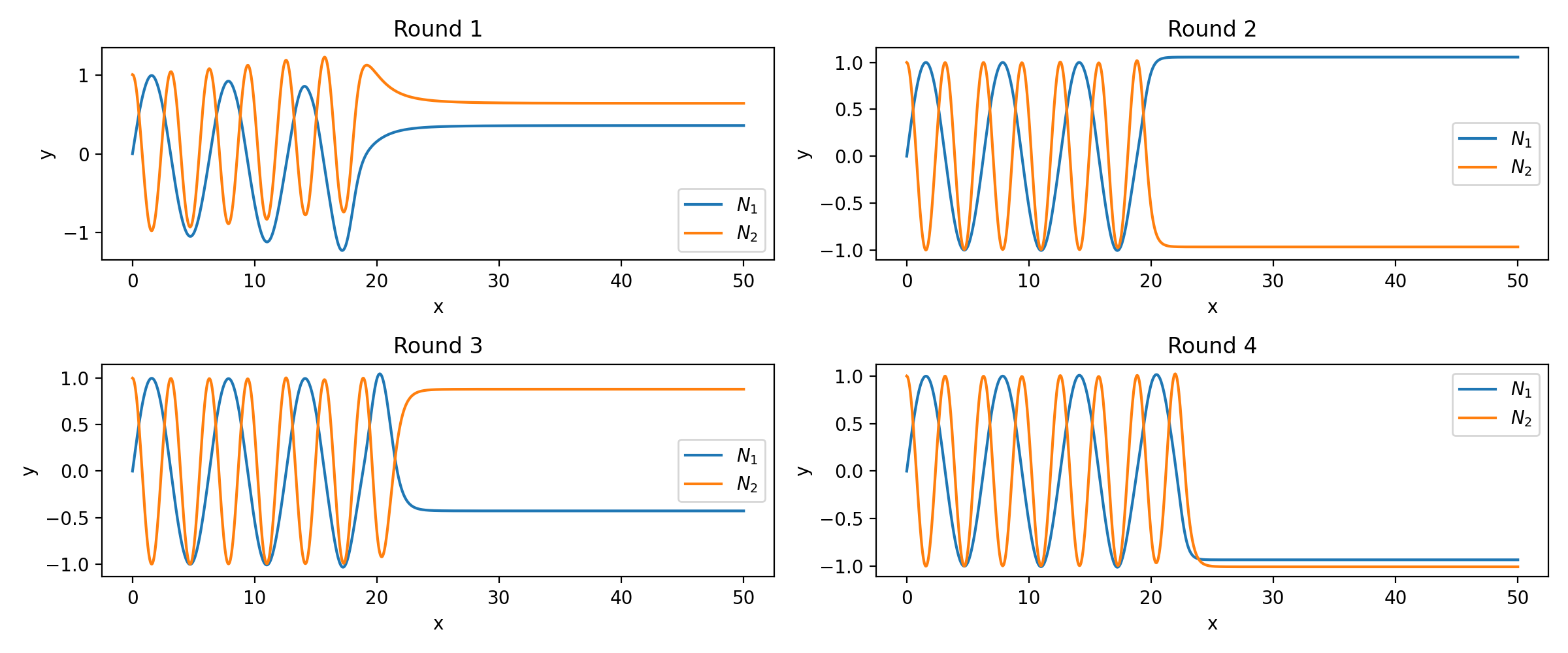}
   \caption{\centering Results of PINN for IVP (\ref{eq:easy}) which is failure to obtain large interval solutions.}
   \label{fig:pinn_easy}
\end{figure}

Consider IVP of an ODEs
\begin{equation}
    \label{eq:easy}
    \begin{cases}
        \frac{\mathrm{d}y_1}{\mathrm{d}x}=\cos x, \\
        \frac{\mathrm{d}y_2}{\mathrm{d}x}=-2\sin 2x,\\
        y_(0)=0,\,\, y_2(0)=1.
    \end{cases}
\end{equation}
It is easy to obtain the analytic solutions of the IVP as $y_1=\sin x, y_2=\cos 2x$ for all $x\in R$. Although this example is very simple, it is not easy for PINN to solve this problem over a large interval due to the periodic oscillation of the solutions. We performed four rounds of pre-training for PINN on the problem, with 10000 iterations per round. The training results are shown in Fig.\ref{fig:pinn_easy}. It can be observed that despite executing four rounds of training, there is still a considerable discrepancy between the training results and the analytic solutions, especially in the interval $[25, 50]$, where the network did not effectively learn the solutions.

\begin{figure}[htbp]
    \centering
    \includegraphics[width=0.9\textwidth]{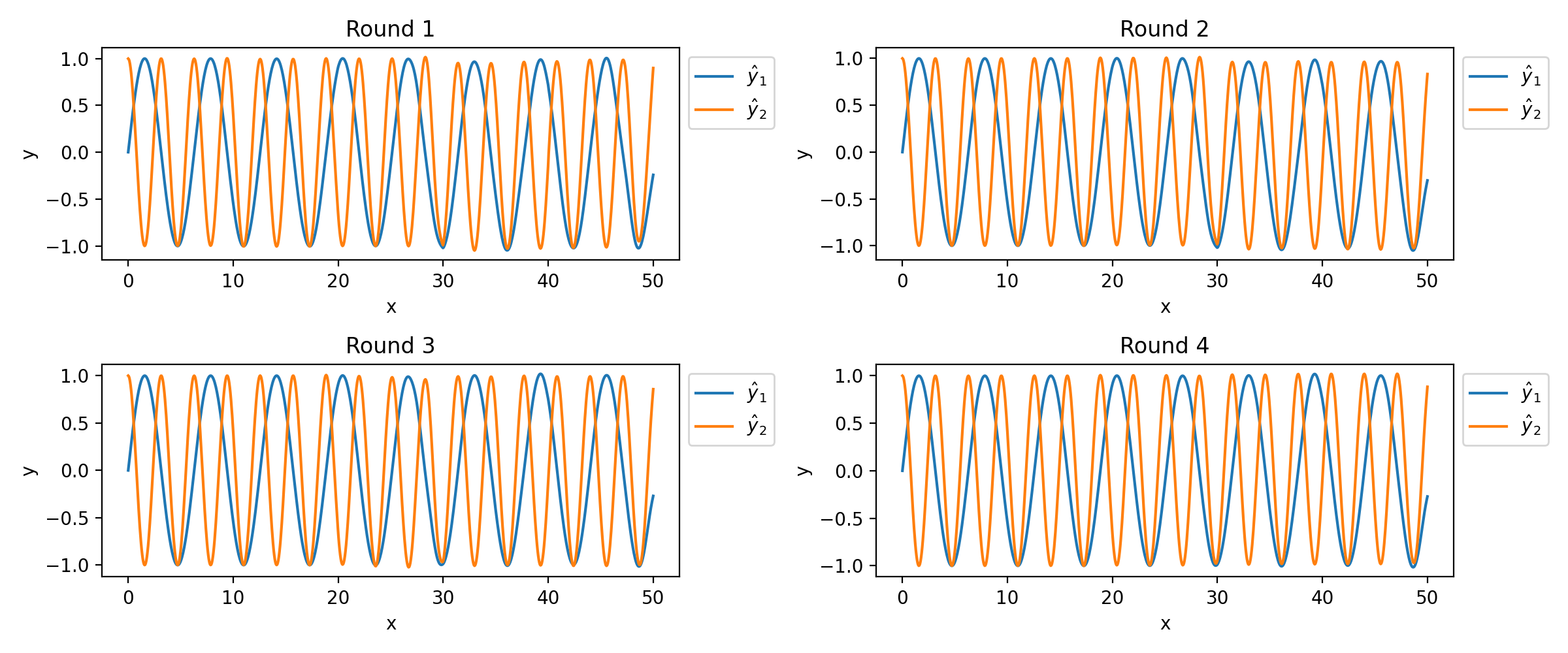}
    \caption{\centering PWNN solutions to problem (\ref{eq:easy}) over [0, 50]}
   \label{fig:pnn_easy}
\end{figure}

By PWNN method, we divide solution interval $[0, 50]$ into five equal sub-intervals and construct the corresponding PWNNs. These networks were also trained in four rounds of pre-training with same iteration number per sub-interval in each training round. The training results are shown in Fig.\ref{fig:pnn_easy}. It can be seen that the PWNNs approximate the analytic solutions quite well after the first training round.

In the experiments, we use uniformly sampled 100 training points in each interval. Table \ref{tab:pnn_pinn} provides detailed comparisons of convergence trends of PINN and PWNN in terms of decreasing of the loss function values during the four training rounds, under same calculation environment. The loss value for the PWNNs refers to the average of the loss function values in each interval. As can be seen from the table, PWNN has obvious advantages in convergence compared with PINN. In Table \ref{tab:pnn_easy}, we present the convergent trend of PWNN as increasing the number of pre-trining rounds, which shows that the convergence of PWNN is improved by significantly decreasing of loss values. 






\begin{table}[htbp]
    \centering
    \caption{Comparison the convergent trends of PWNN and PINN in terms of loss function values.}
    \label{tab:pnn_pinn}
    \begin{tabularx}{\textwidth}{CCC}
    \toprule
    Rounds & PINN & PWNN \\
    \midrule
      1 & 1.59 & 7.88e-3 \\
      2 & 1.49 & 5.14e-3 \\
      3 & 1.42 & 2.93e-3 \\
      4 & 1.34 & 1.92e-3 \\
    \bottomrule
    \end{tabularx}
\end{table}


\begin{table}[htbp]
\centering
\caption{Convergence of PWNNs in each training round.}
\label{tab:pnn_easy}
\begin{tabular}{c|ccccc|c}
\Xhline{1.pt}
\diagbox{Rounds}{loss}{PWNN} & $m_1$ & $m_2$ & $m_3$ & $m_4$ & $m_5$ & Mean value of loss \\
\hline
1 & 1.19e-4 & 6.32e-4 & 1.68e-2 & 3.25e-3 & 1.86e-2 & 7.88e-3\\
2 & 1.68e-5 & 3.42e-4 & 1.57e-2 & 1.48e-3 & 8.18e-3 & 5.14e-3\\
3 & 9.85e-6 & 2.78e-4 & 5.66e-3 & 1.08e-3 & 7.63e-3 & 2.93e-3\\
4 & 6.75e-6 & 2.53e-4 & 8.23e-2 & 9.55e-4 & 7.65e-3 & 1.92e-3\\
\Xhline{1.pt}
\end{tabular}

\end{table}

\subsection{Example 4}
\begin{figure}[htbp]
   \centering
   \includegraphics[width=0.9\textwidth]{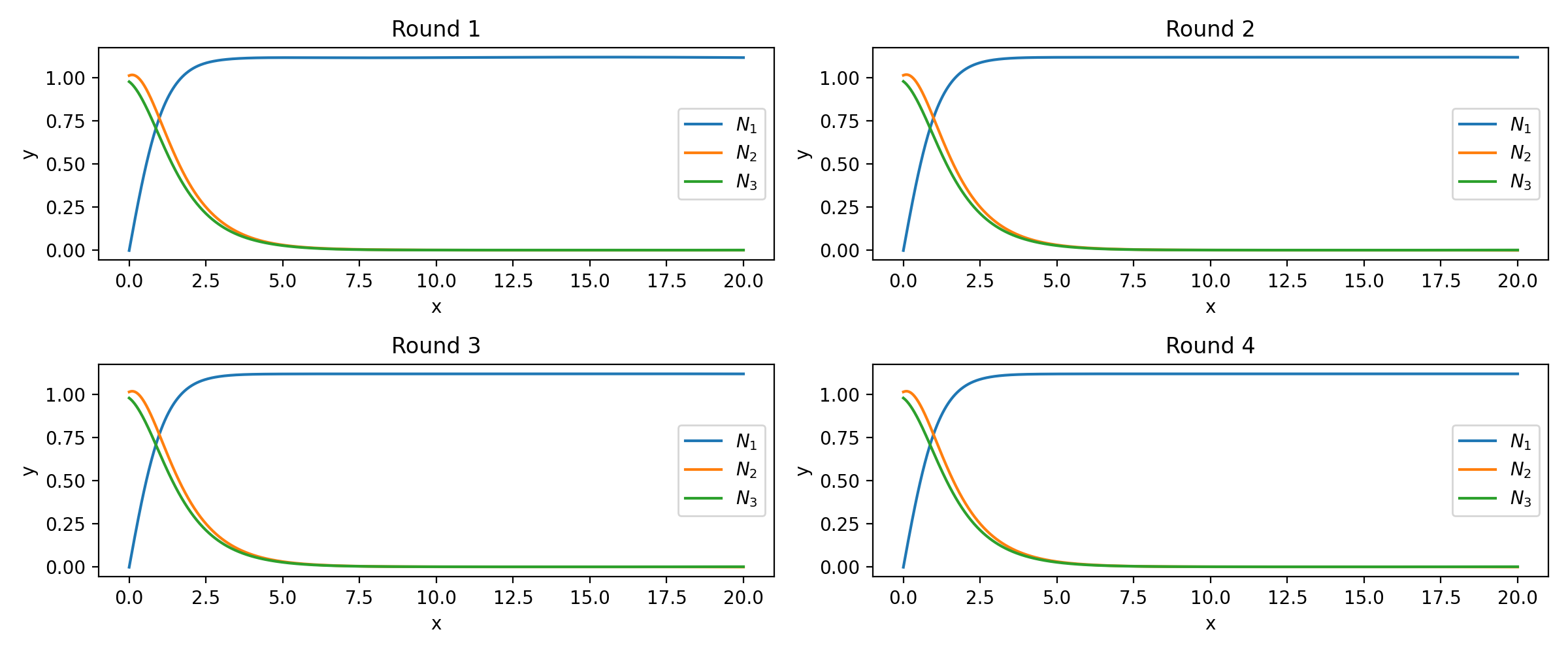}
   \caption{\centering PINN solutions of problem (\ref{eq:e4}) in different rounds}
    \label{fig:pinn_e2}
\end{figure}

Consider the following IVP
\begin{equation}
    \label{eq:e4}
    \begin{cases}
        \frac{\mathrm{d}y_1}{\mathrm{d}x}=y_2y_3,\\
        \frac{\mathrm{d}y_2}{\mathrm{d}x}=-y_1y_3,\\
        \frac{\mathrm{d}y_3}{\mathrm{d}x}=-0.51y_1y_2,\\
        y_1(0)=0,y_2(0)=y_3(0)=1.
    \end{cases}
\end{equation}
As previous examples, we first use PINN to find the numerical solutions of the problem in the interval $[0, 20]$. We perform four rounds of training for PINN, with 10000 iterations per round. The learning rate is set to 0.001 for the first round and 0.0001 for the subsequent rounds. In these experiments, the values of the loss function remains around $1.99*10^{-3}$ for different rounds. This indicates that further increasing the number of training rounds and iterations does not significantly reduce the loss function value. The training results are shown in Fig. \ref{fig:pinn_e2}.

\begin{figure}[htbp]
    \centering
    \includegraphics[width=0.9\textwidth]{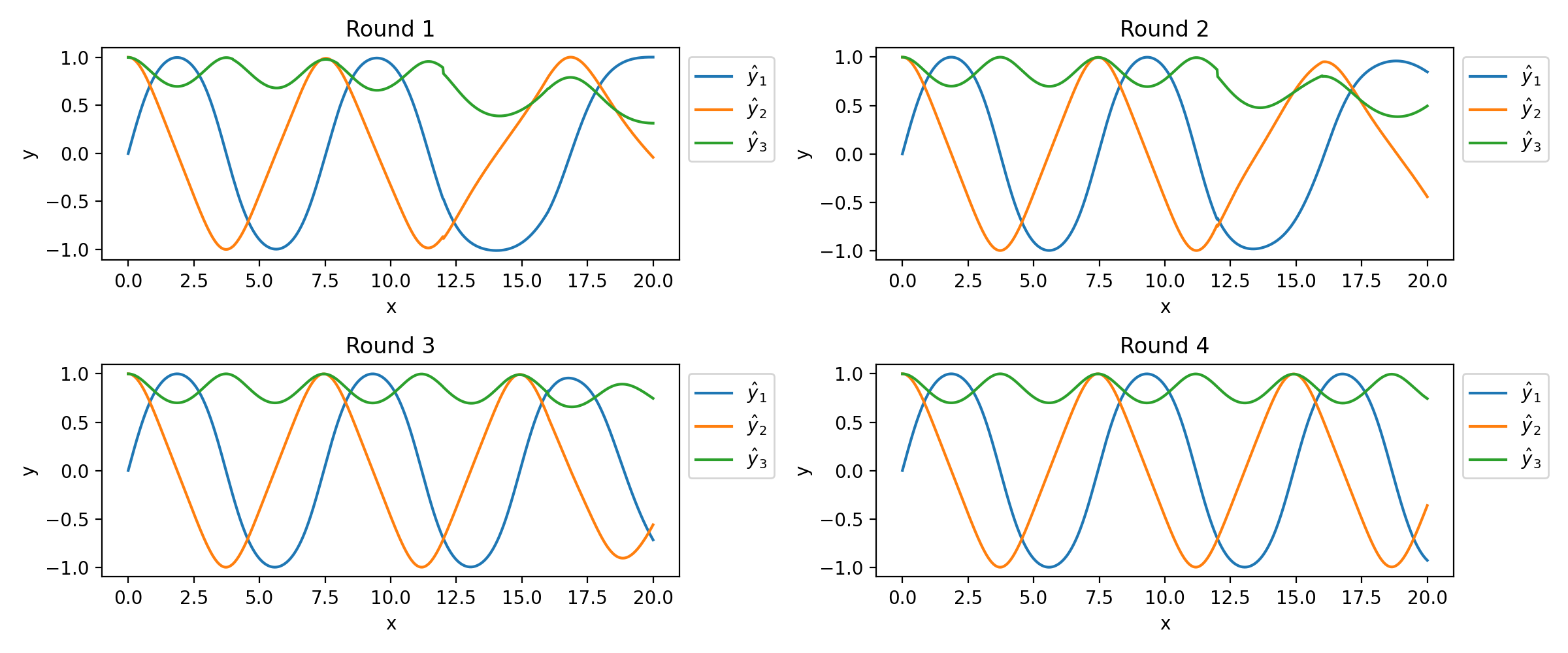}
    \caption{\centering PWNN solutions of IVP (\ref{eq:e4}) in different rounds}
    \label{fig:pnn_e2}
\end{figure}

\begin{figure}[htbp]
    \centering
    \includegraphics[width=0.9\textwidth]{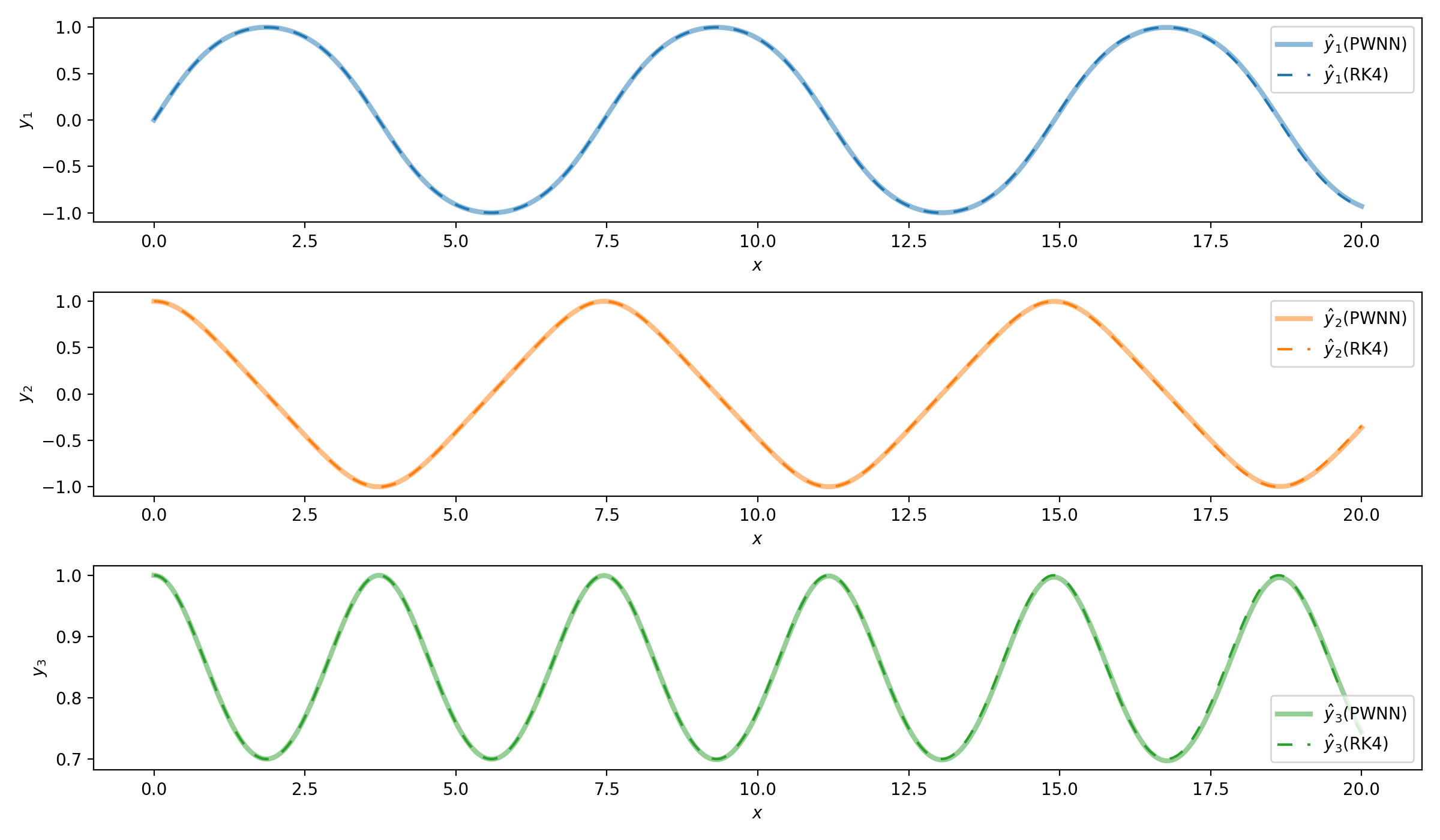}
    \caption{\centering Comparison of solving results of PWNN and RK4 methods for IVP (\ref{eq:e4})}
    \label{fig:pnn_rk4}
\end{figure}

In using PWNN, we divide the solution interval $[0, 20]$ into five equal parts and construct corresponding PWNNs. Then, train each PWNN for four rounds, with each round consisting of 10,000 iterations on each sub-interval. The learning rate for the first round of training is set to $0.001$, while for the remaining rounds, it is set to $0.0001$. The training results of the PWNNs are shown in Fig. \ref{fig:pnn_e2}.
The values of the loss function corresponding to each round of PWNN are shown in Table \ref{tab:pnn_e2}. It can be observed that as the number of training rounds increases, the loss function values of the PWNNs gradually decreases. In Fig. \ref{fig:pnn_rk4}, we show the comparisons of results of the PWNN in the fourth round with the those of RK4. It can be seen that the results of the PWNN method are in good agreement with those of the RK4 method.

\begin{table}[htbp]
\centering
\caption{Comparisons the results obtained by PWNN and RK4 for solving IVP (\ref{eq:e4}).}
\label{tab:pnn_e2}
\begin{tabular}{c|ccccc}
\Xhline{1.pt}
\diagbox{Rounds}{loss}{PWNN} & $m_1$ & $m_2$ & $m_3$ & $m_4$ & $m_5$ \\
\hline
1 & 1.41e-5 & 2.17e-4 & 2.38e-4 & 3.74e-3 & 1.30e-4\\
2 & 1.53e-6 & 1.67e-5 & 1.50e-5 & 5.28e-3 & 1.73e-4 \\
3 & 6.98e-7 & 4.97e-6 & 5.68e-6 & 6.44e-5 & 1.71e-4 \\
4 & 4.43e-7 & 2.20e-6 & 4.09e-6 & 1.94e-5 & 3.61e-5 \\
\Xhline{1.pt}
\end{tabular}

\end{table}

\section{Discussion and Conclusion}
\label{sec:s6}
The neural network-based method for solving differential equations provides solutions with excellent generalization properties. Models based on neural networks offer opportunities to theoretically and practically tackle with differential equation problems across various sciences and engineering applications. While most other techniques offer discrete solutions or solutions with limited differentiability, or strongly depend on the discrete scheme of the variable domain. Effectively using the advantages of nonlinear approximation of artificial neural network and overcoming the shortcomings of traditional numerical methods is a hot research direction at present.

In this paper, based on the advantages of PINN local strong convergence, we propose a piecewise neural network method for solving initial value problems of ODEs over large intervals of the independent variable. This method not only provides an approach for coordination of multiple networks to solve an IVP of ODEs but also introduces a parameter transfer and multiple rounds of pre-training technique to effectively enhance the accuracy of network solutions in training multi-correlation ANN models. On one hand, we offer a new method to address the solution extension problem of initial value problems of differential equations. On the other hand, our work aims to contribute to both the theoretical and practical aspects of applying ANN, providing deep insights into the application of modern methods to traditional mathematical problems.

Through comparative experiments, we prove that under almost the same network training environment, the training time of piecewise neural network is shorter, the convergence speed to the solution of the studied problem is faster, and it can also make up for the defects of the basic ANN algorithm. For example, an interesting observation is that when PINN is used to solve the system of equations (\ref{eq:easy}), the network seems to achieve better training performance on the interval $[0,20]$, while it seems to struggle on the interval $[25,50]$. This phenomenon reflects PINN's tendency to solve the IVP locally rather than the entire interval. Furthermore, in section \ref{sec:4.2}, when PINN is used to solve the equation (\ref{eq:sir}) over a large solution domain, although the final loss value obtained by PINN training is very low, the approximate solution obtained does not match the actual solution. This may reflect two reasons. On the one hand, it shows that the optimization of the residual loss function is only a necessary condition for the artificial neural network output to be an approximate solution of the problem, but not a sufficient condition. On the other hand, this is due to small numerical instabilities during back-propagation because of the complexity of the loss hypersurface where the ANN can settle on a local minimum with a small value for the loss function \cite{Maria}. While the piecewise neural network method proposed in this paper solves these problems effectively to some extend.

Certainly, there are some limitations to the study. We did not consider applying the proposed method to partial differential equation problems, which would involve partitioning high-dimensional domains of independent variables. The successful application of this method strongly depends on the effective training of sub-neural networks on each corresponding interval. In other words, if an ANN in a certain link cannot be effectively trained in the training of a PWNN, the whole algorithm may not complete the solving task. Furthermore, the continuity and differentiability of network solution $\hat{y}(x)$ over the entire interval have yet to be guaranteed theoretically. These are topics for further study in the future.


\bibliography{sn-bibliography}


\begin{thebibliography}{39}
\ifx \bisbn   \undefined \def \bisbn  #1{ISBN #1}\fi
\ifx \binits  \undefined \def \binits#1{#1}\fi
\ifx \bauthor  \undefined \def \bauthor#1{#1}\fi
\ifx \batitle  \undefined \def \batitle#1{#1}\fi
\ifx \bjtitle  \undefined \def \bjtitle#1{#1}\fi
\ifx \bvolume  \undefined \def \bvolume#1{\textbf{#1}}\fi
\ifx \byear  \undefined \def \byear#1{#1}\fi
\ifx \bissue  \undefined \def \bissue#1{#1}\fi
\ifx \bfpage  \undefined \def \bfpage#1{#1}\fi
\ifx \blpage  \undefined \def \blpage #1{#1}\fi
\ifx \burl  \undefined \def \burl#1{\textsf{#1}}\fi
\ifx \doiurl  \undefined \def \doiurl#1{\url{https://doi.org/#1}}\fi
\ifx \betal  \undefined \def \betal{\textit{et al.}}\fi
\ifx \binstitute  \undefined \def \binstitute#1{#1}\fi
\ifx \binstitutionaled  \undefined \def \binstitutionaled#1{#1}\fi
\ifx \bctitle  \undefined \def \bctitle#1{#1}\fi
\ifx \beditor  \undefined \def \beditor#1{#1}\fi
\ifx \bpublisher  \undefined \def \bpublisher#1{#1}\fi
\ifx \bbtitle  \undefined \def \bbtitle#1{#1}\fi
\ifx \bedition  \undefined \def \bedition#1{#1}\fi
\ifx \bseriesno  \undefined \def \bseriesno#1{#1}\fi
\ifx \blocation  \undefined \def \blocation#1{#1}\fi
\ifx \bsertitle  \undefined \def \bsertitle#1{#1}\fi
\ifx \bsnm \undefined \def \bsnm#1{#1}\fi
\ifx \bsuffix \undefined \def \bsuffix#1{#1}\fi
\ifx \bparticle \undefined \def \bparticle#1{#1}\fi
\ifx \barticle \undefined \def \barticle#1{#1}\fi
\bibcommenthead
\ifx \bconfdate \undefined \def \bconfdate #1{#1}\fi
\ifx \botherref \undefined \def \botherref #1{#1}\fi
\ifx \url \undefined \def \url#1{\textsf{#1}}\fi
\ifx \bchapter \undefined \def \bchapter#1{#1}\fi
\ifx \bbook \undefined \def \bbook#1{#1}\fi
\ifx \bcomment \undefined \def \bcomment#1{#1}\fi
\ifx \oauthor \undefined \def \oauthor#1{#1}\fi
\ifx \citeauthoryear \undefined \def \citeauthoryear#1{#1}\fi
\ifx \endbibitem  \undefined \def \endbibitem {}\fi
\ifx \bconflocation  \undefined \def \bconflocation#1{#1}\fi
\ifx \arxivurl  \undefined \def \arxivurl#1{\textsf{#1}}\fi
\csname PreBibitemsHook\endcsname

\bibitem[\protect\citeauthoryear{Braun}{1983}]{ODE}
\begin{bbook}
\bauthor{\bsnm{Braun}, \binits{M.}}:
\bbtitle{Differential {Equations} and {Their} {Applications}}.
\bpublisher{Springer},
\blocation{New York Berlin Heidelbeg Tokyo}
(\byear{1983}).
\doiurl{10.1007/978-1-4684-0173-8}
\end{bbook}
\endbibitem

\bibitem[\protect\citeauthoryear{Piscopo et~al.}{2019}]{Maria}
\begin{barticle}
\bauthor{\bsnm{Piscopo}, \binits{M.L.}},
\bauthor{\bsnm{Spannowsky}, \binits{M.}},
\bauthor{\bsnm{Waite}, \binits{P.}}:
\batitle{Solving differential equations with neural networks: Applications to the calculation of cosmological phase transitions}.
\bjtitle{Phys. Rev. D}
\bvolume{100},
\bfpage{016002}
(\byear{2019})
\doiurl{10.1103/PhysRevD.100.016002}
\end{barticle}
\endbibitem

\bibitem[\protect\citeauthoryear{Yadav et~al.}{2015}]{NUrODE}
\begin{bbook}
\bauthor{\bsnm{Yadav}, \binits{N.}},
\bauthor{\bsnm{Yadav}, \binits{A.}},
\bauthor{\bsnm{Kumar}, \binits{M.}}:
\bbtitle{An {Introduction} to {Neural} {Network} {Methods} for {Differential} {Equations}}.
\bsertitle{{Springer Briefs} in {Applied} {Sciences} and {Technology}}.
\bpublisher{Springer},
\blocation{Springer Dordrecht Heidelberg New York London}
(\byear{2015}).
\doiurl{10.1007/978-94-017-9816-7}
\end{bbook}
\endbibitem

\bibitem[\protect\citeauthoryear{Hornik et~al.}{1989}]{hornik_multilayer_1989}
\begin{barticle}
\bauthor{\bsnm{Hornik}, \binits{K.}},
\bauthor{\bsnm{Stinchcombe}, \binits{M.}},
\bauthor{\bsnm{White}, \binits{H.}}:
\batitle{Multilayer feedforward networks are universal approximators}.
\bjtitle{Neural Networks}
\bvolume{2}(\bissue{5}),
\bfpage{359}--\blpage{366}
(\byear{1989})
\doiurl{10.1016/0893-6080(89)90020-8}
\end{barticle}
\endbibitem

\bibitem[\protect\citeauthoryear{Hornik et~al.}{1990}]{z49}
\begin{barticle}
\bauthor{\bsnm{Hornik}, \binits{K.}},
\bauthor{\bsnm{Stinchcombe}, \binits{M.}},
\bauthor{\bsnm{White}, \binits{H.}}:
\batitle{Universal approximation of an unknown mapping and its derivatives using multilayer feedforward networks}.
\bjtitle{Neural Networks}
\bvolume{3}(\bissue{5}),
\bfpage{551}--\blpage{560}
(\byear{1990})
\doiurl{10.1016/0893-6080(90)90005-6}
\end{barticle}
\endbibitem

\bibitem[\protect\citeauthoryear{Li}{1996}]{z50}
\begin{barticle}
\bauthor{\bsnm{Li}, \binits{X.}}:
\batitle{Simultaneous approximations of multivariate functions and their derivatives by neural networks with one hidden layer}.
\bjtitle{Neurocomputing}
\bvolume{12}(\bissue{4}),
\bfpage{327}--\blpage{343}
(\byear{1996})
\doiurl{10.1016/0925-2312(95)00070-4}
\end{barticle}
\endbibitem

\bibitem[\protect\citeauthoryear{Lagaris et~al.}{1998}]{lagaris_artificial_1998}
\begin{barticle}
\bauthor{\bsnm{Lagaris}, \binits{I.E.}},
\bauthor{\bsnm{Likas}, \binits{A.}},
\bauthor{\bsnm{Fotiadis}, \binits{D.I.}}:
\batitle{Artificial neural networks for solving ordinary and partial differential equations}.
\bjtitle{IEEE Trans. Neural Netw.}
\bvolume{9}(\bissue{5}),
\bfpage{987}--\blpage{1000}
(\byear{1998})
\doiurl{10.1109/72.712178}
\end{barticle}
\endbibitem

\bibitem[\protect\citeauthoryear{Paszke et~al.}{2019}]{paszke_pytorch_2019}
\begin{botherref}
\oauthor{\bsnm{Paszke}, \binits{A.}},
\oauthor{\bsnm{Gross}, \binits{S.}},
\oauthor{\bsnm{Massa}, \binits{F.}},
\oauthor{\bsnm{Lerer}, \binits{A.}},
\oauthor{\bsnm{Bradbury}, \binits{J.}},
\oauthor{\bsnm{Chanan}, \binits{G.}},
\oauthor{\bsnm{Killeen}, \binits{T.}},
\oauthor{\bsnm{Lin}, \binits{Z.}},
\oauthor{\bsnm{Gimelshein}, \binits{N.}},
\oauthor{\bsnm{Antiga}, \binits{L.}},
\oauthor{\bsnm{Desmaison}, \binits{A.}},
\oauthor{\bsnm{K{\"{o}}pf}, \binits{A.}},
\oauthor{\bsnm{Yang}, \binits{E.Z.}},
\oauthor{\bsnm{DeVito}, \binits{Z.}},
\oauthor{\bsnm{Raison}, \binits{M.}},
\oauthor{\bsnm{Tejani}, \binits{A.}},
\oauthor{\bsnm{Chilamkurthy}, \binits{S.}},
\oauthor{\bsnm{Steiner}, \binits{B.}},
\oauthor{\bsnm{Fang}, \binits{L.}},
\oauthor{\bsnm{Bai}, \binits{J.}},
\oauthor{\bsnm{Chintala}, \binits{S.}}:
Pytorch: An imperative style, high-performance deep learning library.
CoRR
\textbf{abs/1912.01703}
(2019)
\doiurl{10.48550/arXiv.1912.01703}
{\href{https://arxiv.org/abs/1912.01703}{{1912.01703}}}
\end{botherref}
\endbibitem

\bibitem[\protect\citeauthoryear{Rall}{1981}]{z46}
\begin{bbook}
\bauthor{\bsnm{Rall}, \binits{L.B.}}:
\bbtitle{Automatic Differentiation: Techniques and Applications}.
\bpublisher{Springer},
\blocation{New York Berlin Heidelbeg}
(\byear{1981}).
\doiurl{10.1007/3-540-10861-0}
\end{bbook}
\endbibitem

\bibitem[\protect\citeauthoryear{Verma}{2000}]{z47}
\begin{barticle}
\bauthor{\bsnm{Verma}, \binits{A.}}:
\batitle{An introduction to automatic differentiation}.
\bjtitle{Current Science}
\bvolume{78}(\bissue{7}),
\bfpage{804}--\blpage{807}
(\byear{2000})
\end{barticle}
\endbibitem

\bibitem[\protect\citeauthoryear{Baydin et~al.}{2018}]{z48}
\begin{barticle}
\bauthor{\bsnm{Baydin}, \binits{A.G.}},
\bauthor{\bsnm{Pearlmutter}, \binits{B.A.}},
\bauthor{\bsnm{Radul}, \binits{A.A.}},
\bauthor{\bsnm{Siskind}, \binits{J.M.}}:
\batitle{Automatic differentiation in machine learning: a survey}.
\bjtitle{Journal of Machine Learning Research}
\bvolume{18}(\bissue{153}),
\bfpage{1}--\blpage{43}
(\byear{2018})
\end{barticle}
\endbibitem

\bibitem[\protect\citeauthoryear{Sirignano and Spiliopoulos}{2018}]{sirignano_dgm_2018}
\begin{barticle}
\bauthor{\bsnm{Sirignano}, \binits{J.}},
\bauthor{\bsnm{Spiliopoulos}, \binits{K.}}:
\batitle{Dgm: A deep learning algorithm for solving partial differential equations}.
\bjtitle{Journal of Computational Physics}
\bvolume{375},
\bfpage{1339}--\blpage{1364}
(\byear{2018})
\doiurl{10.1016/j.jcp.2018.08.029}
\end{barticle}
\endbibitem

\bibitem[\protect\citeauthoryear{Anitescu et~al.}{2019}]{anitescu_artificial_2019}
\begin{barticle}
\bauthor{\bsnm{Anitescu}, \binits{C.}},
\bauthor{\bsnm{Atroshchenko}, \binits{E.}},
\bauthor{\bsnm{Alajlan}, \binits{N.}},
\bauthor{\bsnm{Rabczuk}, \binits{T.}}:
\batitle{Artificial {Neural} {Network} {Methods} for the {Solution} of {Second} {Order} {Boundary} {Value} {Problems}}.
\bjtitle{Computers, Materials \& Continua}
\bvolume{59}(\bissue{1}),
\bfpage{345}--\blpage{359}
(\byear{2019})
\doiurl{10.32604/cmc.2019.06641}
\end{barticle}
\endbibitem

\bibitem[\protect\citeauthoryear{Raissi et~al.}{2019}]{raissi_physics-informed_2019}
\begin{barticle}
\bauthor{\bsnm{Raissi}, \binits{M.}},
\bauthor{\bsnm{Perdikaris}, \binits{P.}},
\bauthor{\bsnm{Karniadakis}, \binits{G.E.}}:
\batitle{Physics-informed neural networks: {A} deep learning framework for solving forward and inverse problems involving nonlinear partial differential equations}.
\bjtitle{Journal of Computational Physics}
\bvolume{378},
\bfpage{686}--\blpage{707}
(\byear{2019})
\doiurl{10.1016/j.jcp.2018.10.045}
\end{barticle}
\endbibitem

\bibitem[\protect\citeauthoryear{Yuan et~al.}{2022}]{yuan_-pinn_2022}
\begin{barticle}
\bauthor{\bsnm{Yuan}, \binits{L.}},
\bauthor{\bsnm{Ni}, \binits{Y.-Q.}},
\bauthor{\bsnm{Deng}, \binits{X.-Y.}},
\bauthor{\bsnm{Hao}, \binits{S.}}:
\batitle{A-{PINN}: {Auxiliary} physics informed neural networks for forward and inverse problems of nonlinear integro-differential equations}.
\bjtitle{Journal of Computational Physics}
\bvolume{462},
\bfpage{111260}
(\byear{2022})
\doiurl{10.1016/j.jcp.2022.111260}
\end{barticle}
\endbibitem

\bibitem[\protect\citeauthoryear{Chiu et~al.}{2022}]{chiu_can-pinn_2022}
\begin{barticle}
\bauthor{\bsnm{Chiu}, \binits{P.-H.}},
\bauthor{\bsnm{Wong}, \binits{J.C.}},
\bauthor{\bsnm{Ooi}, \binits{C.}},
\bauthor{\bsnm{Dao}, \binits{M.H.}},
\bauthor{\bsnm{Ong}, \binits{Y.-S.}}:
\batitle{{CAN}-{PINN}: {A} fast physics-informed neural network based on coupled-automatic–numerical differentiation method}.
\bjtitle{Computer Methods in Applied Mechanics and Engineering}
\bvolume{395},
\bfpage{114909}
(\byear{2022})
\doiurl{10.1016/j.cma.2022.114909}
\end{barticle}
\endbibitem

\bibitem[\protect\citeauthoryear{Huang et~al.}{2022}]{huang_hompinns_2022}
\begin{barticle}
\bauthor{\bsnm{Huang}, \binits{Y.}},
\bauthor{\bsnm{Hao}, \binits{W.}},
\bauthor{\bsnm{Lin}, \binits{G.}}:
\batitle{{HomPINNs}: {Homotopy} physics-informed neural networks for learning multiple solutions of nonlinear elliptic differential equations}.
\bjtitle{Computers $\&$ Mathematics with Applications}
\bvolume{121},
\bfpage{62}--\blpage{73}
(\byear{2022})
\doiurl{10.1016/j.camwa.2022.07.002}
\end{barticle}
\endbibitem

\bibitem[\protect\citeauthoryear{Fang et~al.}{2021}]{fang_data-driven_2021}
\begin{barticle}
\bauthor{\bsnm{Fang}, \binits{Y.}},
\bauthor{\bsnm{Wu}, \binits{G.-Z.}},
\bauthor{\bsnm{Wang}, \binits{Y.-Y.}},
\bauthor{\bsnm{Dai}, \binits{C.-Q.}}:
\batitle{Data-driven femtosecond optical soliton excitations and parameters discovery of the high-order {NLSE} using the {PINN}}.
\bjtitle{Nonlinear Dynamics}
\bvolume{105}(\bissue{1}),
\bfpage{603}--\blpage{616}
(\byear{2021})
\doiurl{10.1007/s11071-021-06550-9}
\end{barticle}
\endbibitem

\bibitem[\protect\citeauthoryear{Bai et~al.}{2022}]{bai_application_2022}
\begin{barticle}
\bauthor{\bsnm{Bai}, \binits{Y.}},
\bauthor{\bsnm{Chaolu}, \binits{T.}},
\bauthor{\bsnm{Bilige}, \binits{S.}}:
\batitle{The application of improved physics-informed neural network ({IPINN}) method in finance}.
\bjtitle{Nonlinear Dynamics}
\bvolume{107}(\bissue{4}),
\bfpage{3655}--\blpage{3667}
(\byear{2022})
\doiurl{10.1007/s11071-021-07146-z}
\end{barticle}
\endbibitem

\bibitem[\protect\citeauthoryear{Meng et~al.}{2020}]{z80}
\begin{barticle}
\bauthor{\bsnm{Meng}, \binits{X.}},
\bauthor{\bsnm{Li}, \binits{Z.}},
\bauthor{\bsnm{Zhang}, \binits{D.}},
\bauthor{\bsnm{Karniadakis}, \binits{G.E.}}:
\batitle{Ppinn: Parareal physics-informed neural network for time-dependent pdes}.
\bjtitle{Computer Methods in Applied Mechanics and Engineering}
\bvolume{370},
\bfpage{113250}
(\byear{2020})
\doiurl{10.1016/j.cma.2020.113250}
\end{barticle}
\endbibitem

\bibitem[\protect\citeauthoryear{Long et~al.}{2019}]{z57}
\begin{barticle}
\bauthor{\bsnm{Long}, \binits{Z.}},
\bauthor{\bsnm{Lu}, \binits{Y.}},
\bauthor{\bsnm{Dong}, \binits{B.}}:
\batitle{Pde-net 2.0: Learning pdes from data with a numeric-symbolic hybrid deep network}.
\bjtitle{Journal of Computational Physics}
\bvolume{399},
\bfpage{108925}
(\byear{2019})
\doiurl{10.1016/j.jcp.2019.108925}
\end{barticle}
\endbibitem

\bibitem[\protect\citeauthoryear{Zha et~al.}{2022}]{z73}
\begin{barticle}
\bauthor{\bsnm{Zha}, \binits{W.}},
\bauthor{\bsnm{Zhang}, \binits{W.}},
\bauthor{\bsnm{Li}, \binits{D.}},
\bauthor{\bsnm{Xing}, \binits{Y.}},
\bauthor{\bsnm{He}, \binits{L.}},
\bauthor{\bsnm{Tan}, \binits{J.}}:
\batitle{Convolution-based model-solving method for three-dimensional, unsteady, partial differential equations}.
\bjtitle{Neural Computation}
\bvolume{34}(\bissue{2}),
\bfpage{518}--\blpage{540}
(\byear{2022})
\doiurl{10.1162/neco_a_01462}
\end{barticle}
\endbibitem

\bibitem[\protect\citeauthoryear{Gao et~al.}{2021}]{z93}
\begin{barticle}
\bauthor{\bsnm{Gao}, \binits{H.}},
\bauthor{\bsnm{Sun}, \binits{L.}},
\bauthor{\bsnm{Wang}, \binits{J.-X.}}:
\batitle{Phygeonet: Physics-informed geometry-adaptive convolutional neural networks for solving parameterized steady-state pdes on irregular domain}.
\bjtitle{Journal of Computational Physics}
\bvolume{428},
\bfpage{110079}
(\byear{2021})
\doiurl{10.1016/j.jcp.2020.110079}
\end{barticle}
\endbibitem

\bibitem[\protect\citeauthoryear{Wang et~al.}{2020}]{z85}
\begin{barticle}
\bauthor{\bsnm{Wang}, \binits{N.}},
\bauthor{\bsnm{Zhang}, \binits{D.}},
\bauthor{\bsnm{Chang}, \binits{H.}},
\bauthor{\bsnm{Li}, \binits{H.}}:
\batitle{Deep learning of subsurface flow via theory-guided neural network}.
\bjtitle{Journal of Hydrology}
\bvolume{584},
\bfpage{124700}
(\byear{2020})
\doiurl{10.1016/j.jhydrol.2020.124700}
\end{barticle}
\endbibitem

\bibitem[\protect\citeauthoryear{Wang et~al.}{2021}]{z86}
\begin{barticle}
\bauthor{\bsnm{Wang}, \binits{N.}},
\bauthor{\bsnm{Chang}, \binits{H.}},
\bauthor{\bsnm{Zhang}, \binits{D.}}:
\batitle{Theory-guided auto-encoder for surrogate construction and inverse modeling}.
\bjtitle{Computer Methods in Applied Mechanics and Engineering}
\bvolume{385},
\bfpage{114037}
(\byear{2021})
\doiurl{10.1016/j.cma.2021.114037}
\end{barticle}
\endbibitem

\bibitem[\protect\citeauthoryear{Zhang and Bilige}{2019}]{zhang_bilinear_2019}
\begin{barticle}
\bauthor{\bsnm{Zhang}, \binits{R.-F.}},
\bauthor{\bsnm{Bilige}, \binits{S.}}:
\batitle{Bilinear neural network method to obtain the exact analytical solutions of nonlinear partial differential equations and its application to p-{gBKP} equation}.
\bjtitle{Nonlinear Dynamics}
\bvolume{95}(\bissue{4}),
\bfpage{3041}--\blpage{3048}
(\byear{2019})
\doiurl{10.1007/s11071-018-04739-z}
\end{barticle}
\endbibitem

\bibitem[\protect\citeauthoryear{Zhang et~al.}{2021}]{zhang2}
\begin{barticle}
\bauthor{\bsnm{Zhang}, \binits{R.}},
\bauthor{\bsnm{Bilige}, \binits{S.}},
\bauthor{\bsnm{Chaolu}, \binits{T.}}:
\batitle{Fractal {Solitons}, {Arbitrary} {Function} {Solutions}, {Exact} {Periodic} {Wave} and {Breathers} for a {Nonlinear} {Partial} {Differential} {Equation} by {Using} {Bilinear} {Neural} {Network} {Method}}.
\bjtitle{Journal of Systems Science and Complexity}
\bvolume{34}(\bissue{1}),
\bfpage{122}--\blpage{139}
(\byear{2021})
\doiurl{10.1007/s11424-020-9392-5}
\end{barticle}
\endbibitem

\bibitem[\protect\citeauthoryear{{Chien-Cheng Yu} et~al.}{2002}]{chien-cheng_yu_adaptive_2002}
\begin{bchapter}
\bauthor{\bsnm{{Chien-Cheng Yu}}},
\bauthor{\bsnm{{Yun-Ching Tang}}},
\bauthor{\bsnm{{Bin-Da Liu}}}:
\bctitle{An adaptive activation function for multilayer feedforward neural networks}.
In: \bbtitle{2002 {IEEE} {Region} 10 {Conference} on {Computers}, {Communications}, {Control} and {Power} {Engineering}. {TENCOM} '02. {Proceedings}.},
pp. \bfpage{645}--\blpage{650}.
\bpublisher{IEEE},
\blocation{Beijing, China}
(\byear{2002}).
\doiurl{10.1109/TENCON.2002.1181357}
\end{bchapter}
\endbibitem

\bibitem[\protect\citeauthoryear{Dushkoff and Ptucha}{2016}]{dushkoff_adaptive_2016}
\begin{barticle}
\bauthor{\bsnm{Dushkoff}, \binits{M.}},
\bauthor{\bsnm{Ptucha}, \binits{R.}}:
\batitle{Adaptive {Activation} {Functions} for {Deep} {Networks}}.
\bjtitle{Electronic Imaging}
\bvolume{28}(\bissue{19}),
\bfpage{1}--\blpage{5}
(\byear{2016})
\doiurl{10.2352/ISSN.2470-1173.2016.19.COIMG-149}
\end{barticle}
\endbibitem

\bibitem[\protect\citeauthoryear{Qian et~al.}{2018}]{qian_adaptive_2018}
\begin{barticle}
\bauthor{\bsnm{Qian}, \binits{S.}},
\bauthor{\bsnm{Liu}, \binits{H.}},
\bauthor{\bsnm{Liu}, \binits{C.}},
\bauthor{\bsnm{Wu}, \binits{S.}},
\bauthor{\bsnm{Wong}, \binits{H.S.}}:
\batitle{Adaptive activation functions in convolutional neural networks}.
\bjtitle{Neurocomputing}
\bvolume{272},
\bfpage{204}--\blpage{212}
(\byear{2018})
\doiurl{10.1016/j.neucom.2017.06.070}
\end{barticle}
\endbibitem

\bibitem[\protect\citeauthoryear{Wang et~al.}{2022}]{wang_when_2022}
\begin{barticle}
\bauthor{\bsnm{Wang}, \binits{S.}},
\bauthor{\bsnm{Yu}, \binits{X.}},
\bauthor{\bsnm{Perdikaris}, \binits{P.}}:
\batitle{When and why {PINNs} fail to train: {A} neural tangent kernel perspective}.
\bjtitle{Journal of Computational Physics}
\bvolume{449},
\bfpage{110768}
(\byear{2022})
\doiurl{10.1016/j.jcp.2021.110768}
\end{barticle}
\endbibitem

\bibitem[\protect\citeauthoryear{Xiang et~al.}{2022}]{xiang_self-adaptive_2022}
\begin{barticle}
\bauthor{\bsnm{Xiang}, \binits{Z.}},
\bauthor{\bsnm{Peng}, \binits{W.}},
\bauthor{\bsnm{Liu}, \binits{X.}},
\bauthor{\bsnm{Yao}, \binits{W.}}:
\batitle{Self-adaptive loss balanced {Physics}-informed neural networks}.
\bjtitle{Neurocomputing}
\bvolume{496},
\bfpage{11}--\blpage{34}
(\byear{2022})
\doiurl{10.1016/j.neucom.2022.05.015}
\end{barticle}
\endbibitem

\bibitem[\protect\citeauthoryear{Jagtap et~al.}{2020}]{jagtap_adaptive_2020}
\begin{barticle}
\bauthor{\bsnm{Jagtap}, \binits{A.D.}},
\bauthor{\bsnm{Kawaguchi}, \binits{K.}},
\bauthor{\bsnm{Karniadakis}, \binits{G.E.}}:
\batitle{Adaptive activation functions accelerate convergence in deep and physics-informed neural networks}.
\bjtitle{Journal of Computational Physics}
\bvolume{404},
\bfpage{109136}
(\byear{2020})
\doiurl{10.1016/j.jcp.2019.109136}
\end{barticle}
\endbibitem

\bibitem[\protect\citeauthoryear{Wen and Chaolu}{2023a}]{Wen1}
\begin{botherref}
\oauthor{\bsnm{Wen}, \binits{Y.}},
\oauthor{\bsnm{Chaolu}, \binits{T.}}:
Learning the nonlinear solitary wave solution of the korteweg–de vries equation with novel neural network algorithm.
Entropy
\textbf{25}(5)
(2023)
\doiurl{10.3390/e25050704}
\end{botherref}
\endbibitem

\bibitem[\protect\citeauthoryear{Wen and Chaolu}{2023b}]{wen2}
\begin{botherref}
\oauthor{\bsnm{Wen}, \binits{Y.}},
\oauthor{\bsnm{Chaolu}, \binits{T.}}:
Study of burgers–huxley equation using neural network method.
Axioms
\textbf{12}(5)
(2023)
\doiurl{10.3390/axioms12050429}
\end{botherref}
\endbibitem

\bibitem[\protect\citeauthoryear{Wen et~al.}{2022}]{wen3}
\begin{barticle}
\bauthor{\bsnm{Wen}, \binits{Y.}},
\bauthor{\bsnm{Chaolu}, \binits{T.}},
\bauthor{\bsnm{Wang}, \binits{X.}}:
\batitle{Solving the initial value problem of ordinary differential equations by lie group based neural network method}.
\bjtitle{PloS one}
\bvolume{17}(\bissue{4}),
\bfpage{0265992}
(\byear{2022})
\doiurl{10.1371/journal.pone.0265992}
\end{barticle}
\endbibitem

\bibitem[\protect\citeauthoryear{Wen and Chaolu}{2023}]{wen4}
\begin{barticle}
\bauthor{\bsnm{Wen}, \binits{Y.}},
\bauthor{\bsnm{Chaolu}, \binits{T.}}:
\batitle{Lie group-based neural networks for nonlinear dynamics}.
\bjtitle{International Journal of Bifurcation and Chaos}
\bvolume{33}(\bissue{14}),
\bfpage{2350161}--\blpage{12}
(\byear{2023})
\doiurl{10.1142/S0218127423501614}
\end{barticle}
\endbibitem

\bibitem[\protect\citeauthoryear{Glorot and Bengio}{2010}]{glorot_xavier_2010}
\begin{barticle}
\bauthor{\bsnm{Glorot}, \binits{X.}},
\bauthor{\bsnm{Bengio}, \binits{Y.}}:
\batitle{Understanding the difficulty of training deep feedforward neural networks}.
\bjtitle{Journal of Machine Learning Research - Proceedings Track}
\bvolume{9},
\bfpage{249}--\blpage{256}
(\byear{2010})
\end{barticle}
\endbibitem

\bibitem[\protect\citeauthoryear{Kingma and Ba}{2017}]{kingma_adam_2017}
\begin{botherref}
\oauthor{\bsnm{Kingma}, \binits{D.P.}},
\oauthor{\bsnm{Ba}, \binits{J.}}:
Adam: {A} {Method} for {Stochastic} {Optimization}.
arXiv.
arXiv:1412.6980 [cs]
(2017).
\doiurl{10.48550/arXiv.1412.6980}
\end{botherref}
\endbibitem

\end{thebibliography}

\end{document}